\documentclass[12pt]{amsart}

\usepackage[utf8]{inputenc}

\usepackage{amsfonts}
\usepackage{amsmath}
\usepackage{amssymb}
\usepackage{mathtools}
\usepackage{tabularx}

\usepackage{geometry} 
\usepackage{pxfonts}
\usepackage{euscript}
\usepackage{bbold,bbm}
\usepackage{comment}
\usepackage{scalerel}
\usepackage{xcolor}
\usepackage{tikz-cd}

\usepackage{hyperref}
\usepackage{graphics}
\usepackage{epstopdf} 
\usepackage[all,2cell,arrow,matrix,tips]{xy} \UseAllTwocells \SilentMatrices
\usepackage{graphicx}
\usepackage{verbatim}
\usepackage{leftidx}
\usepackage{enumerate}
\usepackage{phonetic}
\usepackage{lipsum}
\usepackage{subfigure}
\usepackage{float}

\usepackage[backend=bibtex,style=numeric, maxbibnames=99, minbibnames=99]{biblatex}
\numberwithin{equation}{section}

\theoremstyle{definition}

\newtheorem{thm}{Theorem}[section]
\newtheorem{prop}[thm]{Proposition}
\newtheorem{cor}[thm]{Corollary}
\newtheorem{lem}[thm]{Lemma}
\newtheorem{rem}[thm]{Remark}
\newtheorem{example}[thm]{Example}

\newcommand\Zb {\mathbb{Z}}
\newcommand\Rb {\mathbb{R}}
\newcommand\Cb {\mathbb{C}}

\newcommand\CH {\EuScript{H}}

\newcommand\CK {\EuScript{K}}

\newcommand\CP {\EuScript{P}}

\newcommand{\Cla} {\mathcal{A}}
\newcommand{\Clb} {\mathcal{B}}
\newcommand{\Clc} {\mathcal{C}}

\newcommand{\Clj} {\mathcal{J}}

\newcommand{\Clm} {\mathcal{M}}
\newcommand{\Cln} {\mathcal{N}}

\newcommand{\Clp} {\mathcal{P}}

\newcommand{\Clu} {\mathcal{U}}
\newcommand{\Clv} {\mathcal{V}}
\newcommand{\Clw} {\mathcal{W}}

\newcommand{\Clz} {\mathcal{Z}}

\DeclareMathOperator{\SOT}{SOT}
\DeclareMathOperator{\sot}{sot}
\DeclareMathOperator{\Wot}{WOT}
\DeclareMathOperator{\wot}{wot}
\DeclareMathOperator{\diag}{diag}

\DeclareMathOperator{\Tr}{Tr}

\DeclareMathOperator{\Span}{Span}

\DeclareMathOperator{\supp}{supp}
\DeclareMathOperator{\sa}{sa}
\DeclareMathOperator{\sgn}{sgn}
\DeclareMathOperator{\op}{op}
\DeclarePairedDelimiterX\braketBase[2]{\langle}{\rangle}{%
  #1\,\delimsize\vert\,\mathopen{}#2%
}
\NewDocumentCommand{\braket}{s o m m}{%
  \IfBooleanTF{#1}%
    {\braketBase{#3}{#4}}
    {\IfValueTF{#2}%
      {\braketBase[#2]{#3}{#4}}
      {\braketBase*{#3}{#4}}
    }%
}

\begin{document}

\title{$\Rb$-linear unitary preservers of C$^*$-algebras}

\author{Lu Cui}
\address{College of Science, China University of Petroleum (East China), Qingdao, 266580, China.}
\email{nilu@upc.edu.cn}

\author{Wenming Wu}
\address{School of Mathematical Sciences, Chongqing Normal University, Chongqing, 401331, China.}
\email{wuwm@amss.ac.cn}

\author{Wei Yuan}
\address{Institute of Mathematics, Academy of Mathematics and Systems Science, Chinese Academy of Sciences, Beijing, 100190, China; School of Mathematical Sciences, University of Chinese Academy of Sciences, Beijing 100049, China.}
\email{wyuan@math.ac.cn}
\thanks{Lu Cui is supported by the Natural Science Foundation of Shandong Province under grant number ZR2024QA097. 
Wenming Wu is supported by the NSFC under grant number 12571130 and the Natural Science Foundation of
Chongqing under grant number CSTB2025NSCQLZX0058.
Wei Yuan is supported by the NSFC under grant numbers 12471124 and 12571130.}

\date{}

\begin{abstract}
We study the structure of general $\Rb$-linear unitary preservers between unital C$^*$-algebras. Representing any unital $\Rb$-linear unitary preserver in terms of three $\Rb$-linear maps on the self-adjoint part of the domain yields a canonical decomposition into a Jordan-type summand and a purely non-Jordan summand. When the domain is abelian, we completely classify the purely non-Jordan summand via type~$\mathrm{I}_2$ JW-algebras. For unital C$^*$-algebras with no two-dimensional irreducible representations, we characterize a natural class of purely non-Jordan unitary preservers governed by pairs of Jordan $*$-homomorphisms.
\end{abstract}


\maketitle

\section{Introduction}
The study of linear preserver problems on substructures of operator algebras has a rich history. In general, these problems seek to characterize linear maps that preserve specific algebraic, topological, or geometric substructures or subsets \cite{LiPierce2001}. Among such subsets, the unitary group plays a prominent role, serving as the geometric and algebraic core of operator algebras.

An early result of Marcus \cite{marcus59} states that every $\Cb$-linear, unital map on $M_n(\Cb)$ preserving the unitary group is either a $*$-automorphism or a $*$-antiautomorphism. Explicitly, such a map takes the form
\begin{align*}
 A \mapsto U^*AU \quad \mbox{or} \quad A \mapsto U^*A^tU,
\end{align*} 
where $U$ is a fixed unitary matrix. This result was subsequently generalized to C$^*$-algebras. Russo and Dye ~\cite{RD66} proved that any unital, $\Cb$-linear unitary preserver is a Jordan $*$-homomorphism. Thus $\Cb$-linear unitary preservers essentially determine the Jordan structure of the underlying $C^*$-algebras.

The situation changes substantially when complex linearity is replaced by real linearity. A similar classification for $\Rb$-linear maps preserving the real orthogonal group was obtained in \cite{Wei1975}. More precisely, for $n \neq 2,4,8$, every $\Rb$-linear map on $M_n(\Rb)$ that preserves the orthogonal group $O(n)$ is a composition of maps of the following two types:
\begin{itemize}
\item $A \mapsto U A V$, where $U, V \in O(n)$;
\item $A \mapsto A^{t}$.
\end{itemize}
For $n = 2,4,8$, one must also include an additional family of maps in these compositions, namely those of the form 
\begin{align*}
 \sum_{i,j} a_{ij} E_{ij} \mapsto \sum_{j} a_{1j} U_j,   
\end{align*}
where $E_{ij}$ are the standard matrix units, and ${U_1,\ldots,U_n}$ are pairwise skew-symmetric orthogonal matrices, i.e., $U_i^t U_j + U_j^t U_i = 0$ for all $i \neq j$. These exotic maps indicate that spin systems and type I$_2$ Jordan structures play an essential role in the $\Rb$-linear unitary preserver problem.

The purpose of this paper is to further develop the structural theory of (unital) $\Rb$-linear unitary preservers between unital C$^*$-algebras. Besides the intrinsic interest of this problem, another key motivation stems from the study of the geometry of projection spaces (e.g., \cite{Dye1955,GeherSemrl2016,GeherSemrl2018,GuWuYuan2020,Mori2019,QianWangWuYuan2019, QianShenShiWuYuan2023, ShiShenGuMa2022}), where metric and geometric properties are known to reflect the Jordan structure of the underlying operator algebras. More specifically, let $\Cla$ be a finite factor equipped with its unique normalized trace $\tau$, and consider the Grassmannian of trace $\frac{1}{2}$ projections 
\begin{align*}
  \CP_{\frac{1}{2}}(\Cla)= \left \{P\in\Cla: P=P^*=P^2,\ \tau(P)= \frac{1}{2}  \right \}.
\end{align*}
Unlike Grassmann spaces of projections with trace distinct from $\frac{1}{2}$, the space $\CP_{\frac{1}{2}}$ is distinguished by the fact that the orthocomplementation map $P\mapsto I-P$ leaves it invariant. Moreover, every $P \in \CP_{\frac{1}{2}}(\Cla)$ determines a symmetry $S_P = 2P-I$. As demonstrated by the identities
\begin{align*}
  [S_P,S_Q] = 4[P,Q], \qquad \tau(S_P S_Q) = 4 \tau(PQ)-1,
\end{align*}
both commutativity and transition-probability data for $\CP_{\frac{1}{2}}(\Cla)$ are encoded by this family of unitaries. Therefore, preserver-problems on the unitary group and Grassmannians of projections are naturally linked. In fact, for $0<p<\infty$, every surjective $L^p$-isometry on $\CP_{\frac{1}{2}}(\Cla)$ induces a surjective unital $L^p$-isometry on the unitary group of a compression $P\Cla P$ for some projection $P$ of trace $\frac{1}{2}$; this isometry admits an $\Rb$-linear extension to $P\Cla P$~\cite{QianShenShiWuYuan2023}. This observation makes the structural theory of $\Rb$-linear unitary preservers a linchpin in classifying $L^p$-isometries on Grassmannians in von Neumann algebras. In this work, we lift the surjectivity restriction and study general $\Rb$-linear unitary preservers between unital C$^*$-algebras. For general operator algebra background, we refer the reader to \cite{KR97,T02}.

We now outline the organization of this paper and summarize our main results. Throughout the paper, $\Clc$ denotes a unital C$^*$-algebra, and $\Clb(\CH)$ denotes the algebra of bounded operators on a complex Hilbert space $\CH$. By faithfully representing the codomain, it suffices to consider maps $\varphi:\Clc\to\Clb(\CH)$. We restrict our attention to unital unitary preservers, as any unitary preserver $\varphi$ can be reduced to the unital case by multiplying by $\varphi(I)^*$.

In Section 2, we first show that every unital $\Rb$-linear unitary preserver $\varphi:\Clc\to\Clb(\CH)$ extends uniquely to a normal unital $\Rb$-linear unitary preserver on the enveloping von Neumann algebra $\Clc^{**}$ of $\Clc$ (see Theorem \ref{thm:ext_to_preserves}). This allows us to use von Neumann algebra methods and then restrict the resulting conclusions to $\Clc$. We then establish a structural framework by decomposing such maps in terms of three $\Rb$-linear maps on the self-adjoint part $\Clc_{\sa}$ of $\Clc$. More precisely, every $\Rb$-linear unitary preserver $\varphi: \Clc \to \Clb(\CH)$ admits a representation of the form
\begin{align}\label{equ:varphi_decomp}
  \varphi(H+iK) = \alpha(H)+i\beta_1(H) +iJ\bigl(\alpha(K)+i\beta_2(K)\bigr), \qquad H, K \in \Clc_{\sa},
\end{align}
where $J=-i\varphi(iI)$ is a self-adjoint unitary and the $\Rb$-linear maps $\alpha,\beta_1,\beta_2:\Clc_{\sa}\to\Clb(\CH)_{\sa}$ satisfy a collection of commutativity and multiplicative relations (see Theorem \ref{thm:two_conj_part}). Furthermore, the central projection determined by the range of $\beta_1$ (or $\beta_2$) in the von Neumann algebra generated by the range of $\alpha$ yields a canonical decomposition of $\varphi$ into a Jordan-type summand and a purely non-Jordan summand (see Proposition \ref{prop:Jordan_center_split}). As the Jordan-type summand is determined by Jordan $*$-homomorphisms and the adjoint operation, we focus on the purely non-Jordan summand for the rest of the paper.

In Section 3, we study purely non-Jordan unitary preservers for which the range of $\alpha$ in the decomposition~\eqref{equ:varphi_decomp} is commutative. We show that, for such a map, the image of each map
\begin{align*}
    H \mapsto J\bigl(\alpha(H) + i\beta_j(H)\bigr), \qquad H \in \Clc_{\sa} \quad (j = 1, 2),
\end{align*}
generates a type~$\mathrm{I}_2$ JW-algebra (we refer the reader to~\cite{OS84} for the general theory of JW-algebras). In particular, by decomposing this JW-algebra into type~$\mathrm{I}_{2,\kappa}$ components indexed by cardinals $\kappa$~\cite{Sta82}, we obtain a complete description of $\Rb$-linear unitary preservers when $\Clc$ is abelian (see Theorem \ref{thm:abelian_case}). 

In Section 4, we study unitary preservers on $M_n(\Cb)$ and, unsurprisingly, observe an anomaly for $n=2$. For $n\geq3$, we obtain the well-known classification of unitary preservers on $M_n(\Cb)$ as proved by Marcus \cite{marcus59}, under the weaker condition of only requiring $\Rb$-linearity (see Theorem \ref{thm:Mn_to_Mn}). For $n=2$, exceptional maps appear, which are not of Jordan type. Using the Pauli matrices and their anticommuting relations, we give a classification of $\Rb$-linear unitary preservers from $M_2(\Cb)$ into $\Clb(\CH)$ (see Theorem \ref{thm:classification_M_2}).

In Section 5, we consider unital C$^*$-algebras with no two-dimensional irreducible representations, or equivalently, whose enveloping von Neumann algebras have no type~$\mathrm{I}_2$ direct summand. For any such domain, we completely describe the purely non-Jordan unital $\Rb$-linear unitary preservers satisfying $\beta_1 = \beta_2$ in \eqref{equ:varphi_decomp}. We call these maps binary unitary preservers because they are determined by two Jordan $*$-homomorphisms (see Theorem \ref{thm:cls_binary}). Our proof relies on characterizing pairs of Jordan $*$-homomorphisms on such C$^*$-algebras under certain commutativity relations (see Theorem \ref{thm:decmp_general}). Although this characterization is of independent interest, we defer its proof to the appendix to streamline the exposition.

\section{Structure of $\Rb$-linear unital unitary preservers}

By replacing an $\Rb$-linear unitary preserver $\varphi$ between unital $C^*$-algebras with $\varphi(I)^* \varphi$, the study of $\Rb$-linear unitary preservers reduces to the unital case; hence, without loss of generality, we assume throughout the rest of this paper that all unitary preservers are unital. Moreover, since our results are independent of the specific structure of the codomain, we may apply the GNS construction to assume without loss of generality that the codomain is $\Clb(\CH)$, the von Neumann algebra of bounded linear operators on a Hilbert space $\CH$. In this section, we decompose every such preserver into a sum of $\Rb$-linear maps on self-adjoint operators that satisfy specific algebraic constraints.

We begin by showing that every $\Rb$-linear unitary preserver $\varphi : \Clc \to \Clb(\CH)$ on a unital $C^*$-algebra $\Clc$ extends uniquely to a normal $\Rb$-linear unitary preserver from its enveloping von Neumann algebra $\Clc^{**}$ into $\Clb(\CH)$. This allows us to utilize the abundance of projections in von Neumann algebras hereafter.

Before beginning the proof, we recall some basic facts about mutually skew operators. Let $A, B$ be two operators in a $C^*$-algebra.
$A$ is said to be skew with $B$ if 
\begin{align*}
    A^*B + B^*A = 0, \quad AB^*+BA^* = 0.
\end{align*}

\begin{rem}
It is clear that the skew relation satisfies the following properties:
\begin{itemize}
    \item \textbf{Symmetry:} $A$ is skew with $B$ if and only if $B$ is skew with A.
    Then $A,B$ are said to be \textbf{mutually skew}.
    \item \textbf{Adjoint Invariance:} $A, B$ are mutually skew if and only if $A^*, B^*$ are mutually skew.
    \item If either $A$ or $B$ is a unitary, the two conditions in the definition of the skew relation are equivalent. 
\end{itemize}
\end{rem}

\begin{example}
The Pauli matrices, defined by
\begin{align*}
    \sigma_x = \begin{pmatrix}
       0 & 1\\
       1 & 0 
    \end{pmatrix}, \quad 
    \sigma_y = \begin{pmatrix}
       0 & -i\\
       i & 0 
    \end{pmatrix}, \quad
    \sigma_z = \begin{pmatrix}
       1 & 0\\
       0 & -1 
    \end{pmatrix}
\end{align*} 
are mutually skew.
\end{example}

We restate Lemma 1 from \cite{Wei1975} as follows, and provide a sketch of the proof for the reader's convenience.

\begin{lem}\label{lem:wei2}
    Let $A_1, \ldots, A_n$ be operators in a unital $C^*$-algebra. 
    If the linear combination
    \begin{align*}
        \sum_{j=1}^n a_j A_j 
    \end{align*}
    is unitary for every choice of coefficients $a_j \in \{1, -1\}$, then the operators $A_j$ are mutually skew and 
    \begin{align*}
        \sum_{j=1}^n A_j^*A_j = \sum_{j=1}^n A_j A_j^* = I.
    \end{align*} 
\end{lem}

\begin{proof}
By expanding the identity 
\begin{align*}
   \left ( \sum_{j=1}^n a_j A_j \right )^* \left (\sum_{j=1}^n a_j A_j \right ) = I = \left (\sum_{j=1}^n a_j A_j \right ) \left (\sum_{j=1}^n a_j A_j \right )^*, 
\end{align*}
and varying the values of the coefficients $a_j$, it follows that $A_j$ and $A_k$ are mutually skew whenever $j \neq k$.
\end{proof}

\begin{lem}\label{lem:basic_equ}
Let $\Clc$ be a unital C$^*$-algebra, and let $\varphi: \Clc \to \Clb(\CH)$ be a unital $\Rb$-linear unitary preserver. Then $\|\varphi\| = 1$, and $\varphi$ satisfies 
\begin{align}\label{equ:basic_equ_I}
    \varphi(iH)^* = -\varphi(iH) \quad \text{and} \quad \varphi(iH)^*\varphi(iH) = \frac{1}{2}\bigl(\varphi(H^2) + \varphi(H^2)^*\bigr)
\end{align}
for all $H \in \Clc_{\sa}$.
\end{lem}

\begin{proof}
By \cite[Theorem~1]{kadison85}, $\|\varphi\| = 1$. For each $H \in \Clc_{\sa}$, the Equations \eqref{equ:basic_equ_I} hold because
\begin{align*}
    0 &= \left. \frac{d}{dt} \left( \varphi(e^{itH})^* \varphi(e^{itH}) \right) \right|_{t=0} = \varphi(iH)^* + \varphi(iH), \\
    0 &= \left. \frac{d^2}{dt^2} \left( \varphi(e^{itH})^* \varphi(e^{itH}) \right) \right|_{t=0} = -\varphi(H^2)^* - \varphi(H^2) + 2\varphi(iH)^* \varphi(iH).
\end{align*}   
\end{proof}

In the following, we denote the Jordan product on the real Jordan algebra $\Clc_{\sa}$ of self-adjoint operators in a C$^*$-algebra $\Clc$ by $H \circ K$, defined as
\begin{align*}
    H \circ K = \frac{1}{2}(HK + KH).
\end{align*}

\begin{thm}\label{thm:ext_to_preserves}
Every unital $\Rb$-linear unitary preserver $\varphi : \Clc \to \Clb(\CH)$ admits a unique extension to a normal unital $\Rb$-linear unitary preserver $\widehat{\varphi} : \Clc^{**} \to \Clb(\CH)$.
\end{thm}

\begin{proof}
Note that 
\begin{align*}
   \varphi(A) = \varphi_1(A) + \varphi_2(A^*),
\end{align*}
where 
\begin{align*}
   \varphi_1(A) : = \frac{1}{2}\bigl (\varphi(A) - i\varphi(iA) \bigr), \qquad 
    \varphi_2(A) : = \frac{1}{2}\bigl (\varphi(A^*) + i\varphi(iA^*) \bigr)
\end{align*} 
are bounded $\Cb$-linear maps. Let $\widehat{\varphi}_1$ and $\widehat{\varphi}_2$ be the unique normal extensions of $\varphi_1$ and $\varphi_2$, respectively (for the existence of such extensions, see, for example, Lemma A.2.2 in \cite{BM04}). Then 
\begin{align}
    \widehat{\varphi}(C) := \widehat{\varphi}_1(C) + \widehat{\varphi}_2(C^*), \qquad C \in \Clc^{**}
\end{align}
is the unique unital $\Rb$-linear normal extension of $\varphi$. By Kaplansky's density theorem, we have $\| \widehat{\varphi} \| = 1$.

Note that the map $A \mapsto \varphi(iI)^*\widehat{\varphi}(iA)$ is the $\Rb$-linear normal extension of $C \mapsto \varphi(iI)^* \varphi(iC)$. Let $\rho$ denote either of the unital $\mathbb{R}$-linear normal maps $\widehat{\varphi}$ or $\varphi(iI)^* \widehat{\varphi}(i (\cdot))$. By Lemma \ref{lem:basic_equ}, we have
\begin{align*}
    \rho(iH)^*\rho(iK) + \rho(iK)^*\rho(iH) = \rho(H \circ K) + \rho(H \circ K)^*
\end{align*}
for every $H, K \in \Clc_{\sa}$. Extending the above equation first in $H$ to every operator in $(\Clc^{**})_{\sa}$ with $K$ fixed, and then in $K$ with the other variable fixed, we obtain
\begin{align}\label{equ:ext_to_preserves_I}
   \rho(iL)^* \rho(iL) = \frac{1}{2} \bigl( \rho(L^2) + \rho(L^2)^* \bigr)
\end{align}
for every $L \in (\Clc^{**})_{\sa}$. By Kaplansky's density theorem and Lemma \ref{lem:basic_equ}, we have $\rho(iL)^* = -\rho(iL)$. Combining this with Equation \eqref{equ:ext_to_preserves_I}, we obtain
\begin{align*}
    \|\rho(iL)\xi\|^2 = \|\rho(iL)^*\xi\|^2 = \frac{1}{2} \braket{\xi}{\bigl( \rho(L^2) + \rho(L^2)^* \bigr) \xi}
\end{align*}
for every $\xi \in \CH$. Substituting $\widehat{\varphi}(iL)$ and $-\varphi(iI)^*\widehat{\varphi}(L)$ for $\rho(iL)$, we conclude that the maps
\begin{align*}
   L \mapsto \widehat{\varphi}(iL) \quad \text{and} \quad L \mapsto \widehat{\varphi}(L)
\end{align*}
defined on $(\Clc^{**})_{\sa}$ are both strong$^*$ continuous on bounded sets. This implies that $\widehat{\varphi}$ is strong$^*$ continuous on bounded sets. 

Let $U$ be a unitary in $\Clc^{**}$. By the Borel functional calculus, there exists $L \in (\Clc^{**})_{\sa}$ such that $\|L\| \leq \pi$ and $U = e^{iL}$. By Kaplansky's density theorem, there exists a net $\{L_{\alpha}\} \subset \Clc_{\sa}$ such that $\|L_{\alpha}\| \leq \pi$ and $L_{\alpha} \to L$ strongly. By \cite[Proposition 5.3.2]{KR97}, $e^{iL_{\alpha}} \to e^{iL}= U$ in the strong-$*$ topology. Because multiplication is jointly strongly continuous on bounded sets,
\begin{align*}
    \widehat{\varphi}(U)^* \widehat{\varphi}(U) &= \lim_\alpha \varphi \bigl (e^{iL_\alpha} \bigr )^* \varphi\bigl (e^{iL_\alpha} \bigr ) = I, \\
    \widehat{\varphi}(U) \widehat{\varphi}(U)^* &= \lim_\alpha \varphi \bigl (e^{iL_\alpha} \bigr ) \varphi \bigl (e^{iL_\alpha} \bigr )^* = I.
\end{align*}
Consequently, $\widehat{\varphi}$ is a unitary preserver.
\end{proof}

By Theorem \ref{thm:ext_to_preserves}, the study of $\Rb$-linear unitary preservers on $C^*$-algebras reduces to that on von Neumann algebras, and the resulting conclusions can be restricted back to $C^*$-algebras.

\begin{thm}\label{thm:two_conj_part}
Let $\Clc$ be a unital C$^*$-algebra and let $\varphi: \Clc \to \Clb(\CH)$ be an $\Rb$-linear unital map. Then $\varphi$ is a unitary preserver if and only if $\varphi$ takes the form
\begin{align}\label{equ:two_conj_part_I}
    \varphi(H+iK):= \alpha(H) + i\beta_1(H) + iJ(\alpha(K) + i\beta_2(K)), \quad \forall H,K \in \Clc_{\sa},
\end{align}
where $J \in \Clb(\CH)$ is a self-adjoint unitary, and $\alpha, \beta_1, \beta_2: \Clc_{\sa} \to \Clb(\CH)_{\sa}$ are $\Rb$-linear maps satisfying the following conditions:  
    \begin{enumerate}
        \item[I.] $\alpha(I) = I$;
       \item[II.] For every $H \in \Clc_{\sa}$ and $j \in \{1, 2\}$, we have $[J, \alpha(H)] = 0$ and
        \begin{align*}
            J \circ \beta_j(H) = 0, \qquad \alpha(H^2) &= \alpha(H)^2 + \beta_j(H)^2. 
        \end{align*} 
        (Note that, by polarization, $\alpha(H\circ K) = \alpha(H) \circ \alpha(K) + \beta_j(H) \circ \beta_j(K)$ for every $H, K \in \Clc_{\sa}$.)
        \item[III.] For every pair of commuting operators $H, K \in \Clc_{\sa}$,
            \begin{align*}
            [\alpha(H), \alpha(K)] = [\alpha(H), \beta_1(K)] = [\alpha(H), \beta_2(K)] = [\beta_1(H), \beta_2(K)] = 0;
        \end{align*}  
    \end{enumerate}
\end{thm}

\begin{proof}
   Sufficiency follows by direct calculation, so we establish only necessity. By Theorem \ref{thm:ext_to_preserves}, it suffices to prove the necessity under the assumption that $\Clc$ is a von Neumann algebra. 
   
   Let $J = -i\varphi(iI)$ and
   \begin{align*}
        \alpha(H) &:= \frac{1}{2} \bigl (\varphi(H) + \varphi(H)^* \bigr ), \\
        \beta_1(H) &:= \frac{1}{2i} \bigl (\varphi(H) - \varphi(H)^* \bigr ), \\
        \beta_2(H) &:= -\frac{1}{2} \bigl (J\varphi(iH) - \varphi(iH)J \bigr )
   \end{align*}
   for every $H \in \Clc_{\sa}$. By Lemma \ref{lem:basic_equ}, $J$ is self-adjoint and $\beta_2(H)^* = \beta_2(H)$. Combining the first equation in \eqref{equ:basic_equ_I} with the polarization of the second equation with respect to $I$ and $H$, we obtain
    \begin{align*}
        \alpha(H) = -\frac{i}{2} \bigl( J \varphi(iH) + \varphi(iH)J \bigr).
    \end{align*}
    Therefore
    \begin{align*}
    \varphi(H+ iK) = \alpha(H) + i \beta_1(H) + iJ\bigl (\alpha(K) + i\beta_2(K) \bigr ).    
    \end{align*} 

    It is clear that $\alpha(I) = I$, $[J, \alpha(H)] = 0$ and $J \circ \beta_2(H) = 0$. By Lemma \ref{lem:basic_equ},
    \begin{align*}
        \alpha(H^2) = \alpha(H)^2 + \beta_2(H)^2 + i\bigl[\alpha(H), \beta_2(H)\bigr].
    \end{align*}
    Conjugating this identity by $J$ and using the fact that $J$ commutes with $\alpha(H)$ and $\alpha(H^2)$ while anticommuting with $\beta_2(H)$, we obtain
    \begin{align*}
       [\alpha(H), \beta_2(H)] = 0, \mbox{ and }  \alpha(H^2) = \alpha(H)^2 + \beta_2(H)^2.
    \end{align*}
    Since $-iJ\varphi(i(\cdot))$ admits the same decomposition as $\varphi$ with $\beta_1$ and $\beta_2$ interchanged, the remaining equations in (II) follow by applying the above result to $-iJ\varphi(i(\cdot))$.

    Since self-adjoint operators can be approximated by real linear combinations of projections, to prove (III) it suffices to show that  
    \begin{align}\label{equ:two_conj_part_II}
        [\alpha(P), \alpha(Q)] = [\alpha(P), \beta_1(Q)] = [\alpha(P), \beta_2(Q)] = [\beta_1(P), \beta_2(Q)] = 0
    \end{align}  
    for every pair of commuting projections $P, Q \in \Clc$. Applying Lemma \ref{lem:wei2} to $\cos \theta (2P-I)$ and $i\sin \theta (2Q-I)$, we conclude that $\varphi(2P-I)$ and $\varphi(i(2Q-I))$ are mutually skew. Using the equalities in (II) together with the fact that $\varphi(2P-I)$ and $\varphi(i(2Q-I))$ are mutually skew, we obtain
    \begin{align*}
        \varphi(i(2Q-I)) \varphi(P) \varphi(i(2Q-I))^* = \varphi(P)^*.
    \end{align*}
    This implies that 
    \begin{align*}
        \alpha(P) \bigl(\alpha(Q) + i\beta_2(Q) \bigr) &= \bigl(\alpha(Q) + i\beta_2(Q) \bigr) \alpha(P), \\ 
        \beta_1(P) \bigl(\alpha(Q) + i\beta_2(Q) \bigr) &= \bigl(\alpha(Q) + i\beta_2(Q) \bigr)\beta_1(P).
    \end{align*}
    Conjugating these identities by $J$ yields
    \begin{align*}
        \alpha(P) \bigl(\alpha(Q) - i\beta_2(Q) \bigr) &= \bigl(\alpha(Q) - i\beta_2(Q) \bigr) \alpha(P), \\ 
        \beta_1(P) \bigl(-\alpha(Q) + i\beta_2(Q) \bigr) &= \bigl(-\alpha(Q) + i\beta_2(Q) \bigr)\beta_1(P).
    \end{align*}
    Combining these sets of identities then yields the equalities in \eqref{equ:two_conj_part_II}. 
\end{proof}

\begin{rem}\label{rem:pauli_dilation}
Let $\Clc$ and $\Clb$ be unital $C^*$-algebras, and let $\alpha, \beta_1, \beta_2: \Clc_{\sa} \to \Clb_{\sa}$ be $\Rb$-linear maps satisfying conditions (I) and (III) of Theorem \ref{thm:two_conj_part}, and 
\begin{align*}
    \alpha(H^2) = \alpha(H)^2 + \beta_j(H)^2, \qquad \forall H \in \Clc_{\sa}, j \in \{1, 2\}.
\end{align*}
Then 
\begin{align*}
    \varphi(H+iK) :=&\alpha(H)\otimes I_2+i\beta_1(H)\otimes\sigma_z + i(I \otimes \sigma_x)(\alpha(K) \otimes I_2 + i\beta_2(K) \otimes \sigma_z) 
\end{align*}
is a unital $\Rb$-linear unitary preserver from $\Clc$ to $\Clb\otimes M_2(\Cb)$.   
\end{rem}

For the rest of this section, we assume $\varphi: \Clc \to \Clb(\CH)$ is an $\Rb$-linear unital unitary preserver with the decomposition given in \eqref{equ:two_conj_part_I}. If $\Clc$ is a von Neumann algebra, we denote its projection lattice by $\CP(\Clc)$.

\begin{cor}\label{cor:beta_zero}
Assume that $J = 2F-I$ for some projection $F \in \Clb(\CH)$. If $\beta_1$ and $\beta_2$ are zero maps, then there exist $\Cb$-linear Jordan $*$-homomorphisms $\phi_1: \Clc \to F\Clb(\CH)F$ and $\phi_2: \Clc \to (I-F)\Clb(\CH)(I-F)$ such that
    \begin{align}\label{equ:beta_zero_I}
        \varphi(A) = \phi_1(A) + \phi_2(A^*), \quad \forall A \in \Clc. 
    \end{align} 
\end{cor}

\begin{proof}
    Since $\beta_1$ and $\beta_2$ are zero maps, Theorem \ref{thm:two_conj_part} implies that $\varphi(H)^* = \varphi(H)$ and $\varphi(H^2) = \varphi(H)^2$ for every $H \in \Clc_{\sa}$. The decomposition \eqref{equ:beta_zero_I} then follows by defining $\phi_1(A) = \varphi(A)F$ and $\phi_2(A) = \varphi(A^*)(I-F)$, noting that $\varphi(\cdot)F$ is $\Cb$-linear and $\varphi(\cdot)(I-F)$ is conjugate-linear.
\end{proof}

The $\Rb$-linear unital unitary preserver $\varphi: \Clc \to \Clb(\CH)$ is said to be of \textbf{pure Jordan-type} if it is of the form \eqref{equ:beta_zero_I}. For every projection $P$ in the commutant of the image of $\varphi$ satisfying $P \beta_i(\cdot) \equiv 0$ for $i = 1, 2$, the direct summand $\varphi_P(A) = \varphi(A)P$ is a unitary preserver from $\Clc$ into $P\Clb(\CH)P$ of pure Jordan-type. An $\Rb$-linear unital unitary preserver $\varphi$ is called \textbf{purely non-Jordan} if it possesses no non-zero direct summands of pure Jordan type.

\begin{rem}
If the range of a unital $\Rb$-linear unitary preserver is contained in an abelian C$^*$-algebra, then $\beta_1$ and $\beta_2$ are zero maps since
\begin{align*}
    \alpha(K) + i\beta_2(K) &= J \bigl( J(\alpha(K) + i\beta_2(K)) \bigr)\\ 
    &= \bigl( J(\alpha(K) + i\beta_2(K)) \bigr) J = \alpha(K) - i\beta_2(K).
\end{align*}
In particular, unital $\Rb$-linear unitary preservers with commutative ranges must be of pure Jordan type.
\end{rem}

\begin{cor}\label{cor:surjective_pure_jordan}
If $\varphi(\Clu(\Clc))$ is WOT-dense in $\Clu(\Clb)$ for a von Neumann algebra $\Clb$, then $\varphi$ is of pure Jordan-type.
\end{cor}

\begin{proof}
     Assume that $\varphi(\Clu(\Clc))$ is $\Wot$-dense in $\Clu(\Clb)$ for a von Neumann algebra $\Clb$. Since every element in a von Neumann algebra is an $\Rb$-linear combination of unitaries in itself, the WOT-closure of $\varphi(\Clc)$ is $\Clb$. For $L \in \Clb_{\sa}$, let $H_\lambda + iK_\lambda \in \Clc$ be a net such that $\varphi(H_\lambda + iK_\lambda)$ tends to $L$ in $\Wot$. Since the adjoint operation is WOT-continuous, we have
    \begin{align*}
    JLJ &= \lim_{\lambda} \frac{J(\varphi(H_\lambda +iK_\lambda)+\varphi(H_\lambda + iK_\lambda)^*)J}{2} \\
    &= \lim_{\lambda} J\alpha(H_\lambda)J = \lim_{\lambda} \alpha(H_\lambda) = L.
    \end{align*}
    Therefore $J$ is in the center of $\Clb$, which implies it commutes with $\varphi(\Clc)$. 
    The equation
    \begin{align*}
        \alpha(H) + i\beta_j(H) = J(\alpha(H) + i\beta_j(H))J = \alpha(H) - i\beta_j(H),
    \end{align*}
    implies that $\beta_j(H) = 0$, $j = 1, 2$. Thus, $\varphi$ is of pure Jordan-type by Corollary \ref{cor:beta_zero}. 
\end{proof}

Let $\Clw^*(\alpha)$ be the von Neumann algebras generated by the range of $\alpha$.

\begin{prop}\label{prop:Jordan_center_split}
We have the equality of projections
\begin{equation}\label{equ:Jordan_center_split}
    E:= \bigvee_{H \in \Clc_{\mathrm{sa}}} \supp(\beta_1(H)) = \bigvee_{H \in \Clc_{\mathrm{sa}}} \supp(\beta_2(H)) \in \Clw^*(\alpha).
\end{equation}
Moreover, the central support of $E$ in $\Clw^*(\alpha)$ is the unique projection $Q$ in the commutant of the range of $\varphi$ such that $Q\varphi$ is purely non-Jordan and $(I-Q)\varphi$ is of pure Jordan type. In particular, $JE = EJ$, and $\varphi$ is purely non-Jordan if and only if the central support of $E$ in $\Clw^*(\alpha)$ is $I$. Furthermore, if $\Clc$ is a von Neumann algebra, then 
\begin{align}\label{equ:Jordan_center_split_I}
   E = \bigvee_{P \in \CP(\Clc)} \supp (\beta_j(P)), \qquad j = 1, 2.
\end{align}
\end{prop}

\begin{proof}
    For every $H \in \Clc_{\sa}$, we have 
    \begin{align*}
        \supp (\beta_1(H)) = \supp (\beta_1(H)^2) = \supp (\alpha(H^2) - \alpha(H)^2) = \supp (\beta_2(H)).
    \end{align*}
    Consequently, the projection  
    \begin{align*}
       \bigvee_{H \in\Clc_{\sa}} \supp\left(\beta_1(H)\right) = \bigvee_{H \in\Clc_{\sa}} \supp\left(\beta_2(H)\right)  
    \end{align*}
    lies in $\Clw^*(\alpha)$. And $J$ commutes with $E$ since $J \in \Clw^*(\alpha)'$. Note that $\beta_j$ is norm-continuous. The equalities in \eqref{equ:Jordan_center_split_I} then follow by approximating self-adjoint operators with $\Rb$-linear combinations of projections. 
    
    Let $P$ be the central support of $E$. Then $P$ commutes with $J$ and the ranges of $\beta_j$ and $\alpha$ since $J \in \Clw^*(\alpha)'$ and $(I-P)\beta_j(H)=0$ for $H\in\Clc_{\sa}$ and $j=1,2$. Assume that $Q$ is a projection in $\Clw^*(\alpha)'$ that satisfies $Q \beta_j \equiv 0$. Then 
    \begin{align*}
        Q A\beta_j(H)\xi = AQ \beta_j(H)\xi = 0
    \end{align*}
    for every $A \in \Clw^*(\alpha)$, $H \in \Clc_{\sa}$, and $\xi \in \CH$. Since  
    \begin{align*}
        P \CH = \overline{\Span \{A\beta_j(H)\xi: A \in \Clw^*(\alpha), H \in \Clc_{\sa}, \xi \in \CH\}}, 
    \end{align*}
    we must have $Q \leq I-P$. This implies that $P$ is the unique projection in $\varphi(\Clc)'$ such that $P\varphi$ is purely non-Jordan and $(I-P)\varphi$ is of pure Jordan-type. 
\end{proof}

\begin{rem}\label{rem:ext_non_Jordan}
Let $\widehat{\varphi}: \Clc^{**} \to \Clb(\CH)$ be the normal extension of $\varphi$. For all $H, K \in (\Clc^{**})_{\sa}$, we have
\begin{align*}
    \widehat{\varphi}(H + i K) = \widehat{\alpha}(H) + i\widehat{\beta_1}(H) + iJ\bigl(\widehat{\alpha}(K) + i\widehat{\beta_2}(K)\bigr),
\end{align*}
where 
\begin{align*}
    \widehat{\alpha}(H) &= \frac{1}{2} \bigl(\widehat{\varphi}(H) + \widehat{\varphi}(H)^*\bigr),\\ 
    \widehat{\beta_1}(H) &= \frac{1}{2i} \bigl(\widehat{\varphi}(H) - \widehat{\varphi}(H)^*\bigr),\\
    \widehat{\beta_2}(H) &= -\frac{1}{2} \bigl(J\widehat{\varphi}(iH) - \widehat{\varphi}(iH)J\bigr)
\end{align*}
are the normal extensions of $\alpha$, $\beta_1$, and $\beta_2$, respectively. Note that $\Clw^*(\alpha) = \Clw^*(\widehat{\alpha})$ and 
\begin{align*}
    \bigvee_{H \in (\Clc^{**})_{\sa}} \supp \widehat{\beta_j}(H) = \bigvee_{H \in \Clc_{\sa}} \supp \beta_j(H), \qquad j = 1, 2.
\end{align*}
Consequently, $\widehat{\varphi}$ is purely non-Jordan if and only if $\varphi$ is purely non-Jordan.
\end{rem}

\section{$\Rb$-linear unital unitary preservers on abelian domains}

In this section, we provide a complete characterization of $\Rb$-linear purely non-Jordan unital unitary preservers from an abelian unital $C^*$-algebra into $\Clb(\CH)$. In fact, we consider a slightly more general setting by requiring only that the images of self-adjoint operators have pairwise commuting self-adjoint parts. The following is an example of such a map.

\begin{example}\label{example:cond_exp_const}
    Let $\Cla$ be a separable finite von Neumann algebra equipped with a faithful normal tracial state $\tau$, and $\Clb$ an abelian von Neumann subalgebra of $\Cla$. 
    Consequently, $L^2(\Cla, \tau)$ forms a right $\Clb$-module, and the set of left $\Clb$-bounded vectors in $L^2(\Cla, \tau)$ constitutes an inner-product $\Clb$-module. 
    Let $\{\xi_k\}_{k = 0}^{\infty}$ be a Pimsner-Popa basis for this module with $\xi_0 = I$. 
    
    Assume there exists a family of mutually anti-commuting self-adjoint unitaries $\{J\} \cup \{S_n\}_{n = 1}^{\infty} \subset \Clb' \cap \Cla$. 
    For any $H \in \Cla_{\sa}$, define $\alpha(H) := \braket{I}{H}_{\Clb}$ and
    \begin{align*}
        \beta_1(H) &:= \sum_{k=1}^{\infty} \left( \frac{\braket{\xi_k}{H}_{\Clb} + \braket{\xi_k}{H}^*_{\Clb}}{2} S_{4k-3} + \frac{\braket{\xi_k}{H}_{\Clb} - \braket{\xi_k}{H}^*_{\Clb}}{2i} S_{4k-2} \right), \\ 
        \beta_2(H) &:= -iJ \sum_{k=1}^{\infty} \left( \frac{\braket{\xi_k}{H}_{\Clb} + \braket{\xi_k}{H}^*_{\Clb}}{2} S_{4k-1} + \frac{\braket{\xi_k}{H}_{\Clb} - \braket{\xi_k}{H}^*_{\Clb}}{2i} S_{4k} \right).
    \end{align*}
    Then unitary $J$ and the maps $\alpha, \beta_1$, and $\beta_2$ satisfy the conditions of Theorem \ref{thm:two_conj_part}.
\end{example}

Throughout the rest of this section, $\varphi$ denotes an $\Rb$-linear purely non-Jordan unital unitary preserver from a unital C$^*$-algebra $\Clc$ to $\Clb(\CH)$ with the decomposition given in \eqref{equ:two_conj_part_I}. We further impose the following commutativity condition on $\alpha$:
\begin{align*}
    [\alpha(H), \alpha(K)] = 0, \quad \forall H, K \in \Clc_{\sa}.
\end{align*} 
Let $\Clz$ denote the WOT-closed $\Rb$-linear abelian algebra generated by the image of $\alpha$. In particular, by Proposition \ref{prop:Jordan_center_split}, we have
\begin{align*}
    \bigvee_{H \in \Clc_{\sa}} \supp (\beta_1(H)) = \bigvee_{H \in \Clc_{\sa}} \supp (\beta_2(H)) = I.
\end{align*}

\begin{lem}\label{lem:alpha_beta_commute}
    For every $K \in \Clc_{\sa}$, $\beta_j(K) \in \Clz'$, $j = 1, 2$.
\end{lem}

\begin{proof}
It suffices to show that
\begin{align*}
    [\alpha(H), \beta_j(K)] = 0, \quad \forall H, K \in \Clc_{\sa}. 
\end{align*}
By Theorem \ref{thm:two_conj_part}, $[\alpha(L), \beta_j(L)] = 0$ for all $L \in \Clc_{\sa}$. 
Polarizing this with $L = H+K$ yields
    \begin{align*}
        0 = [\alpha(H+K), \beta_j(H+K)] = [\alpha(H), \beta_j(K)] + [\alpha(K), \beta_j(H)].
    \end{align*}
Therefore
\begin{align*}
    [\alpha(H), [\alpha(H), \beta_j(K)]] = -[\alpha(H), [\alpha(K), \beta_j(H)]] = 0.
\end{align*}
By the Kleinecke-Shirokov theorem (see, for example, \cite{kl57}), $[\alpha(H), \beta_j(K)]$ is quasinilpotent, and since it is skew-adjoint, $[\alpha(H), \beta_j(K)] = 0$.
\end{proof}

For $j \in \{1, 2\}$, let
\begin{align*}
    \Clv_{j, 0} := \Span_{\Rb} \{Z \beta_j(H) : Z \in \Clz,\, H \in \Clc_{\sa}\}, \qquad \Clv_j := \overline{\Clv_{j,0}}^{\wot}.
\end{align*}
By Theorem \ref{thm:two_conj_part} and Lemma \ref{lem:alpha_beta_commute}, we have $K \circ L \in \Clz$ for all $K, L \in \Clv_{j}$.

\begin{lem}\label{lem:moduel_stru}
    For each $j \in \{1, 2\}$, $\Clv_j$ is a $\Clz$-module equipped with the $\Clz$-valued inner product
    \begin{align}\label{eq:module_str_I}
        \braket{A}{B}_{\Clz} := A \circ B, \qquad A, B \in \Clv_j.
    \end{align}
    Furthermore, there exists a unique $\SOT$-continuous isometric $\Clz$-module isomorphism $\Phi\colon \Clv_2 \to \Clv_1$ satisfying
    \begin{align}\label{eq:module_str_II}
        \Phi\bigl(Z \beta_2(H)\bigr) = Z \beta_1(H), \quad Z \in \Clz, H \in \Clc_{\sa}.
    \end{align}
    Moreover, $[A, \Phi(A)] = 0$ for every $A \in \Clv_2$.
\end{lem}

\begin{proof}
    By Lemma \ref{lem:alpha_beta_commute} and Theorem \ref{thm:two_conj_part}, for any $Z_k, Z'_l \in \Clz$ and $H_k, H'_l \in \Clc_{\sa}$, we have
    \begin{align*}
        \braket{\sum_{k} Z_k \beta_j(H_k)}{\sum_{l} Z'_l \beta_j(H'_l)}_{\Clz} 
        &= \sum_{k,l} Z_k Z'_l \big( \beta_j(H_k) \circ \beta_j(H'_l) \big) \\
        &= \sum_{k,l} Z_k Z'_l \big( \alpha(H_k \circ H'_l) - \alpha(H_k) \circ \alpha(H'_l) \big) \in \Clz.
    \end{align*}
    Because $\Clv_j = \overline{\Clv_{j,0}}^{\sot}$, the separate $\SOT$-continuity of operator multiplication implies that \eqref{eq:module_str_I} defines a $\Clz$-valued inner product on $\Clv_j$, $j = 1, 2$.
 
    Let $\Phi_0 : \Clv_{2,0} \to \Clv_{1,0}$ be the $\Rb$-linear map given by $\Phi_0(Z \beta_2(H)) := Z \beta_1(H)$. Because the right-hand side of the above equation is independent of $j$, we obtain
    \begin{align}\label{equ:module_str_III}
        \Phi_0(B) \circ \Phi_0(C) = B \circ C, \qquad \forall B, C \in \Clv_{2,0}.
    \end{align}
    In particular, $\|\Phi_0(B)\xi\| = \|B\xi\|$ for every $\xi \in \CH$. Thus, $\Phi_0$ extends uniquely to an $\SOT$-continuous $\Clz$-module map $\Phi: \Clv_2 \to \Clv_1$ satisfying
    \begin{align*}
        \Phi(A)\xi = \lim_\lambda \Phi(A_\lambda)\xi
    \end{align*}
    for every $\xi \in \CH$ and every net $\{A_\lambda \} \subset \Clv_{2,0}$ converging to $A \in \Clv_2$ in $\SOT$. By \eqref{equ:module_str_III}, $\Phi(A)^2 = A^2$ for every $A \in \Clv_2$, which means that $\Phi$ is an isometric $\Clz$-module isomorphism whose inverse can be defined analogously.
   
    Finally, for every $B = \sum_{k=1}^n Z_k \beta_2(H_k) \in \Clv_{2,0}$, Theorem~\ref{thm:two_conj_part} implies 
    \begin{align*}
        [B, \Phi(B)] &= \sum_{k,l} Z_k Z_l [\beta_2(H_k), \beta_1(H_l)] \\
        &= \frac{1}{2} \sum_{k,l} Z_k Z_l \big( [\beta_2(H_k), \beta_1(H_l)] + [\beta_2(H_l), \beta_1(H_k)] \big) = 0.
    \end{align*}
    By polarization, we obtain $[B, \Phi(C)] + [C, \Phi(B)] = 0$ for every $B, C \in \Clv_{2,0}$. Taking $\SOT$-limits separately in each variable, we conclude that $[A, \Phi(A)] = 0$ for every $A \in \Clv_2$.
\end{proof}

\begin{lem}\label{lem:supp_beta_rel}
There exists a self-adjoint unitary $V \in \Clv_j$, for each $j = 1, 2$.
\end{lem}

\begin{proof}
    Note that  
    \begin{align*}
        I =  \bigvee_{H \in \Clc_{\sa}} \supp (\beta_j(H)) \leq \bigvee_{T \in \Clv_j} \supp(T).   
    \end{align*}
    Let $\Clp = \{\supp(T) : T \in \Clv_j \} \subset \Clz$. By Zorn's Lemma, there exists a maximal family of mutually orthogonal non-zero projections $\{E_k\}_{k \in \Lambda} \subseteq \Clp$. For each $k \in \Lambda$, there exists self-adjoint operator $T_k \in \Clv_j$ such that $\|T_k\| \leq 1$ and $\supp(T_k) = E_k$. Let $T := \sum_{k \in \Lambda} E_k T_k \in \Clv_j$. Then $\supp(T)= \sum_{k \in \Lambda} E_k \in \Clz$. If $\supp(T) \neq I$, there exists $K \in \Clc_{\sa}$ such that the projection $F = \supp((I - \supp(T))\beta_j(K)) \neq 0$. Then $F$ is a non-zero projection in $\Clp$ that is orthogonal to every $E_k$, contradicting the maximality of $\{E_k\}_{k \in \Lambda}$. Let $T = |T|V$ be the polar decomposition of $T$. Then $V$ is a self-adjoint unitary in $\Clv_j$, since
    \begin{align*}
        V = \SOT-\lim_{n \to \infty} \left(T^2 + \frac{1}{n}I\right)^{-1/2} T.      
    \end{align*} 
\end{proof}

Let $\Clj$ be the JW-algebra generated by the image of the $\Rb$-linear map $\Clc_{\sa} \to \Clb(\CH)_{\sa}$ defined by
\begin{align*}
    H \mapsto J(\alpha(H) + i\beta_2(H)), \qquad H \in \Clc_{\sa}.
\end{align*} 
We refer the reader to \cite{OS84} for general background on JW-algebras.

\begin{lem}\label{lem:Jordan_sub_alg}
    The JW-algebra $\Clj$ is the direct sum $\Clj = \Clz \oplus \Clz J \oplus iJ\Clv_2$, and it is of type I$_2$ with center $\Clz$. 
\end{lem}

\begin{proof}
    Note that $J \in \Clj$. By Theorem \ref{thm:two_conj_part}, we have
    \begin{align*}
        \alpha(H) = J \circ J(\alpha(H) + i\beta_2(H)) \in \Clj.
    \end{align*}
    Therefore, $\Clz$, $\Clz J$, $iJ \Clv_2 \subset \Clj$. For every $H, K \in \Clc_{\sa}$ and $W, Z \in \Clz$,  
    \begin{align*}
        &WJ \circ ZJ = WZ \in \Clz;\\
        &(iWJ\beta_2(H)) \circ (iZJ\beta_2(K)) = WZ(\alpha(H \circ K) - \alpha(H)\alpha(K)) \in \Clz;\\ 
        &W \circ ZJ = WZJ \in \Clz J;\\
        &W \circ (iZJ\beta_2(H)) = iWZJ\beta_2(H) \in iJ \Clv_2.
    \end{align*}
    Therefore, the $\Rb$-linear span of $\Clz$, $\Clz J$, and $iJ\Clv_2$ is closed under the Jordan product. We next show that the $\Rb$-linear space is a direct sum of the three WOT-closed subspaces, and thus a JW-algebra.

    By Lemma \ref{lem:supp_beta_rel}, there exists a self-adjoint unitary $V \in \Clv_2$, which satisfies $JVJ = -V$ by Theorem \ref{thm:two_conj_part}. Then, for $A \in \Clz$, $B \in \Clz J$, and $C \in iJ\Clv_2$, the maps
    \begin{align*}
        A + B + C &\mapsto \frac{1}{2}\bigl[(A + B + C) - J(A + B + C)J\bigr] = C, \\
        A + B + C &\mapsto \frac{1}{4}\bigl[(A + B + C) + J(A + B + C)J + V(A + B + C)V + VJ(A + B + C)JV\bigr] = A, \\
        A + B + C &\mapsto \frac{1}{4}\bigl[(A + B + C) + J(A + B + C)J - V(A + B + C)V - VJ(A + B + C)JV\bigr] = B
    \end{align*}
    are WOT-continuous projections onto $iJ\Clv_2$, $\Clz$, and $\Clz J$, respectively, implying that the $\Rb$-linear span of $\Clz$, $\Clz J$, and $iJ\Clv_2$ is a WOT-closed direct sum. Since it is closed under the Jordan product and contains the generators of $\Clj$, we obtain
    \begin{align*}
        \Clj =\Clz \oplus \Clz J \oplus iJ\Clv_2.
    \end{align*}

    Let $A \in \Clz$, $B \in \Clz J$, $C \in iJ\Clv_2$, and let $V \in \Clv_2$ be a self-adjoint unitary. Setting $P = \frac{1}{2}(I + J)$, we have
    \begin{align*}
        P(A + B + C)P = (A + B)P, \qquad (I - P)(A + B + C)(I - P) = (A + B)(I - P).
    \end{align*}
    The commutativity of $\Clz \oplus \Clz J$ implies that $P$ and $I - P$ are abelian projections in $\Clj$. Since $VJ = -JV$, the symmetry $iJV \in \Clj$ satisfies
    \begin{align*}
        (iJV) P (iJV) = I - P,
    \end{align*}
    which shows that $\Clj$ is a JW-algebra of type I$_2$. Finally, assume that $A + B + C$ is in the center of $\Clj$. Since $[J, A + B + C] = 0$, we have $C = 0$. Moreover, $[iJV, A + B] = 0$ implies $B = 0$. Consequently, the center of $\Clj$ is $\Clz$.  
\end{proof}

Recall that a spin system in a JW-algebra $\Clm$ is a family of at least two symmetries (self-adjoint unitaries) $\{S_j\}_{j \in J} \subset \Clm$ that mutually anticommute. A spin system $\{S_j\}_{j \in J}$ is said to be locally maximal if for every non-zero central projection $E \in \Clm$, the family $\{S_j E\}_{j \in J}$ is a maximal spin system in the reduced JW-algebra $\Clm E$. A JW-algebra $\Clm$ of type $\text{I}_2$ is said to be of type I$_{2,\kappa}$ for a cardinal $\kappa$ if it contains a locally maximal spin system $\{S_j\}_{j \in J}$ with cardinality $|J| = \kappa$. In this case, relative to $\{S_j\}_{j \in J}$, every element $L \in \Clm$ admits a unique decomposition
\begin{align*}
    L = L_0 + \sum_{j \in J} L_j S_j,
\end{align*}
where $L_0$ and each $L_j$ belong to the center of $\Clm$, and the sum is understood as the $\SOT$-limit of finite sums. Furthermore, every JW-algebra of type I$_2$ decomposes uniquely as a direct sum of JW-algebras of type I$_{2,\kappa}$ over distinct cardinals $\kappa$ (see \cite[Theorem 6.3.14]{OS84} or \cite[Theorem 2]{Sta82}).

Let 
\begin{align*}
    \Clj = \bigoplus_{\kappa \in \Lambda} \Clj_{\kappa}
\end{align*}
be the decomposition of $\Clj$ into a direct sum of type I$_{2,\kappa}$ JW-algebras, and let $Q_{\kappa} \in \Clz$ be the central projection corresponding to the summand $\Clj_{\kappa}$. In what follows, by identifying the cardinal $\kappa$ with its underlying index set, we define $\kappa_0 = \kappa \setminus \{0\}$.

\begin{lem}\label{lem:non_Jordan_alpha_abelian}
For each $\kappa \in \Lambda$, the $\Rb$-linear maps $Q_{\kappa}\beta_1, Q_{\kappa}\beta_2: \Clc_{\sa} \to \Clb(Q_{\kappa}\CH)$ admit the form
\begin{align*}
    Q_{\kappa}\beta_1(H) = \sum_{j \in \kappa_0} Z_j(H)S_{1, j}, \qquad Q_{\kappa}\beta_2(H) = \sum_{j \in \kappa_0} Z_j(H) S_{2,j},
\end{align*}
where each $Z_j : \Clc_{\sa} \to Q_{\kappa}\Clz$ is an $\Rb$-linear map, and $\{S_{1, j}\}_{j \in \kappa_0} \subseteq Q_{\kappa}\Clv_1$ and $\{S_{2,j}\}_{j \in \kappa_0} \subseteq Q_{\kappa}\Clv_2$ are families of self-adjoint unitaries in $\Clb(Q_{\kappa}\CH)$ satisfying the following conditions:
\begin{enumerate}
    \item $\{J Q_{\kappa}\} \cup \{S_{1,j}\}_{j \in \kappa_0}$ and $\{J Q_{\kappa}\} \cup \{S_{2, j}\}_{j \in \kappa_0}$ are spin systems in $\Clb(Q_{\kappa}\CH)$;
    \item $Q_{\kappa}\alpha(H^2) = (Q_{\kappa}\alpha(H))^2 + \sum_{j \in \kappa_0} Z_{j}(H)^2$ for all $H \in \Clc_{\sa}$;
    \item for every pair of commuting operators $H, K \in \Clc_{\sa}$, 
    \begin{align*}
        \left[\sum_{j \in \kappa_0} Z_j(H)S_{1,j}, \sum_{k \in \kappa_0} Z_{k}(K)S_{2, k} \right] = 0.
    \end{align*}
\end{enumerate}
\end{lem}

\begin{proof}
    By replacing $\varphi$ with $Q_{\kappa} \varphi$, we may assume without loss of generality that $\Clj = \Clj_{\kappa}$ and $Q_{\kappa} = I$, and we omit it in the sequel to simplify notation. 

    By Lemma \ref{lem:Jordan_sub_alg}, there exists a family of self-adjoint unitaries $\{S_{2,j}\}_{j \in \kappa_0} \subseteq \Clv_2$ such that 
    \begin{align*}
        \{J\} \cup \{iJS_{2,j}\}_{j \in \kappa_0}     
    \end{align*} 
    is a locally maximal spin system in $\Clj$. In particular, $\{J\} \cup \{S_{2, j}\}_{j \in \kappa_0}$ is a spin system. Let $Z_j(K) := \beta_2(K) \circ S_{2,j}$. We have
    \begin{align}\label{eq:phi_K_decomp}
        \beta_2(K) = \sum_{j \in \kappa_0} Z_j(K) S_{2,j},
    \end{align}
    since the equality $iJ\beta_2(K) = \sum_{j \in \kappa_0} Z_j(K) (iJS_{2,j})$ holds in $\Clj$. Applying Lemma \ref{lem:moduel_stru}, we obtain a family of pairwise anticommuting self-adjoint unitaries $\{S_{1,j}\}_{j \in \kappa_0} \subseteq \Clv_1$ such that 
    \begin{align*}
        \beta_1(H) = \sum_{j \in \kappa_0} Z_j(H) S_{1,j}.
    \end{align*} 
    Since $J \beta_1(H) J = -\beta_1(H)$, it follows that $\{J\} \cup \{S_{1,j}\}_{j \in \kappa_0}$ forms a spin system. Finally, conditions (2), (3) are ensured by Theorem \ref{thm:two_conj_part}.
\end{proof}

If $\Clc$ is a unital abelian C$^*$-algebra, Theorem \ref{thm:two_conj_part} implies that the range of $\alpha$ is commutative. Consequently, we obtain the following description of $\Rb$-linear purely non-Jordan unital unitary preservers on unital abelian C$^*$-algebras.

\begin{thm}\label{thm:abelian_case}
    Let $\varphi$ be a purely non-Jordan unital $\Rb$-linear unitary preserver from a unital abelian C$^*$-algebra $\Clc$ into $\Clb(\CH)$. The range of the map $-i\varphi(i(\cdot)): \Clc_{\sa} \to \Clb(\CH)$ generates a JW-algebra $\Clj$ of type I$_2$. Suppose that $\Clj$ decomposes as
    \begin{align*}
        \Clj = \bigoplus_{\kappa \in \Lambda} \Clj_{\kappa},
    \end{align*}
    where each $\Clj_{\kappa}$ is a JW-algebra of type I$_{2,\kappa}$. Let $Q_{\kappa}$ be the central projection in the center $\Clz$ of $\Clj$ corresponding to the summand $\Clj_{\kappa}$. Then $Q_{\kappa}\varphi$ takes the form
    \begin{align*}
        Q_{\kappa}\varphi(H+iK) = Z_{0}(H) + i \sum_{j \in \kappa_0} Z_j(H) S_{1,j} + i J_{\kappa} \left( Z_0(K) + i \sum_{j\in \kappa_0} Z_j(K)S_{2,j} \right),
    \end{align*}
    where 
    \begin{align*}
        \{J_{\kappa} \} \cup \{S_{1,j}\}_{j \in \kappa_0} \cup \{S_{2,j}\}_{j \in \kappa_0} 
    \end{align*} 
    is a family of self-adjoint unitaries in $(Q_{\kappa} \Clz)' \cap \Clb(Q_{\kappa}\CH)$, and $\{Z_j\}_{j \in \kappa}$ is a family of $\Rb$-linear maps from $\Clc_{\sa}$ to $Q_{\kappa}\Clz$ satisfying
    \begin{enumerate}
        \item $Z_0(I) = Q_{\kappa}$ and 
        \begin{align*}
           \bigvee_{H \in \Clc_{\sa}} \supp \left ( \sum_{j \in \kappa_0} Z_j(H)^2  \right ) = Q_{\kappa}; 
        \end{align*}
        \item $Z_0(H^2) = Z_0(H)^2 + \sum_{j\in \kappa_0} Z_{j}(H)^2$, for all $H \in \Clc_{\sa}$;
        \item $\{J_{\kappa} \} \cup \{S_{1,j}\}_{j \in \kappa_0} \cup \{iJ_{\kappa} S_{2,j}\}_{j \in \kappa_0}$ is a spin system. 
    \end{enumerate}

\end{thm}

\begin{proof}
    Assume that $\varphi$ admits the decomposition given in \eqref{equ:two_conj_part_I}. Set $J_{\kappa} := Q_{\kappa}J$ and $Z_0 := Q_{\kappa}\alpha$, and let $\{S_{1,j}\}_{j \in \kappa_0}, \{S_{2,j}\}_{j \in \kappa_0}$, and $\{Z_j\}_{j \in \kappa_0}$ be the families of self-adjoint unitaries and $\Rb$-linear maps obtained in Lemma \ref{lem:non_Jordan_alpha_abelian}. It remains only to show that condition (3) holds. By Theorem \ref{thm:two_conj_part}, $[A, B] = 0$ for all $A \in \Clv_1$ and $B \in \Clv_2$. In particular, $[S_{1,j}, S_{2,k}] = 0$ for all $j, k \in \kappa_0$, which verifies condition (3).
\end{proof}

\section{$\Rb$-linear unital unitary preservers on $M_n(\Cb)$}

In this section, we first show that every unital $\Rb$-linear unitary preserver on $M_n(\Cb)$ must be of pure Jordan-type if $n \geq 3$. We then provide a complete description of such maps on $M_2(\Cb)$.

\begin{thm}\label{thm:Mn_to_Mn}
    Let $\varphi: M_n(\Cb) \to M_n(\Cb)$ be a unital $\Rb$-linear unitary preserver. If $n \geq 3$, then there exists a unitary $U \in M_n(\Cb)$ such that $\varphi$ takes exactly one of the following forms:
    \begin{itemize}
        \item $\varphi(A)= U^* A U$, 
        \item $\varphi(A)= U^* A^T U$,  
        \item $\varphi(A)= U^* A^* U$, 
        \item $\varphi(A)= U^* \overline{A}U$. 
    \end{itemize}
\end{thm}

\begin{proof}
Suppose $\varphi$ is a unital $\Rb$-linear unitary preserver with the decomposition given in Theorem \ref{thm:two_conj_part}
Since $\alpha$ is a real-linear commutativity-preserving map on $M_n(\Cb)_{\sa}$, by Theorem 2 in \cite{choi87}, either the image of $\alpha$ is contained in an abelian subalgebra, or $\alpha$ takes the forms
\begin{itemize}
    \item  $\alpha(H) = c U H U^* + f(H)I$,
    \item $\alpha(H) = c U H^t U^* + f(H)I$,
\end{itemize} 
where $U$ is a unitary, $c \in \Rb$, and $f$ is an $\Rb$-linear functional on $M_n(\Cb)_{\sa}$.

We first consider the case $\alpha(H) = c U H U^* + f(H)I$ with $c \neq 0$. Since $[J, \alpha(H)] = 0$ for every $H \in M_n(\Cb)_{\sa}$, we have $[U^*JU, H] = 0$. Thus $J = \pm I$. By Corollary \ref{cor:beta_zero}, $\varphi$ is induced by a Jordan $*$-isomorphism. Applying \cite[Theorem 3.3]{stormer65}, it follows that $\varphi$ takes the form stated in the theorem.

Now, assume the image of $\alpha$ is abelian and $J\neq \pm I$. Since $[J, \alpha(H)] = 0$, we can simultaneously diagonalize the image of $\alpha$ and $J$. Without loss of generality, we assume that 
\begin{align*}
    J = \begin{pmatrix}
        I_k & 0 \\
        0 & -I_{n-k}
    \end{pmatrix},\qquad 1 \leq k \leq n-1.
\end{align*}
Recall that $\alpha$ is unital. We have $\alpha(H) = \diag(\rho_1(H), \dots, \rho_n(H))$, where $\rho_i$ is a state on $M_n(\Cb)$ for $i = 1, \ldots, n$. The anti-commutativity $J \beta_1(H) = -\beta_1(H) J$ implies that 
\begin{align*}
    \beta_1(H) = \begin{pmatrix} 
        0 & \tilde{\beta}_1(H) \\ 
        \tilde{\beta}_1(H)^* & 0 \end{pmatrix},
\end{align*}
where $\tilde{\beta}_1(H) \in M_{k, n-k}(\Cb)$. 
Since
\begin{align*}
    \alpha(H^2) - \alpha(H)^2 = \beta_1(H)^2 = 
     \begin{pmatrix} 
        \tilde{\beta}_1(H)\tilde{\beta}_1(H)^* & 0 \\ 0 & \tilde{\beta}_1(H)^*\tilde{\beta}_1(H) 
    \end{pmatrix},
\end{align*}
we have $\tilde{\beta}_1(H)\tilde{\beta}_1(H)^* = \diag (\rho_1(H^2) - \rho_1(H)^2, \ldots, \rho_k(H^2) - \rho_k(H)^2)$.

Let $r_1(H)$ be the first row of $\tilde{\beta}_1(H)$ (which is a vector in $\Cb^{n-k}$). Note that the rank of the real quadratic form 
\begin{align*}
   H \mapsto  \|r_1(H)\|^2= \rho_1(H^2) - \rho_1(H)^2
\end{align*}
is at least $2(n-1)$. However, since $r_1$ maps $M_n(\Cb)_{\sa}$ into $\Cb^{n-k}$, the rank of the quadratic form $\|r_1(H)\|^2$ is at most $2(n-k)$. 
Therefore, the inequality $2(n-1) \le 2(n-k)$ implies $k \le 1$.


Similarly, by considering the first column of $\tilde{\beta}_1(H)$, a symmetric argument yields $2(n-1) \leq 2k$, which implies $n - 1 \leq k$. Combining these gives $n \leq 2$, contradicting the assumption that $n \geq 3$. Thus, the abelian image case is impossible.
\end{proof}

For the remainder of this section, we provide a complete classification of unital $\Rb$-linear unitary preservers with domain $M_2(\Cb)$. \\

Let $\varphi\colon M_2(\Cb) \to \Clb(\CK)$ be a unital $\Rb$-linear unitary preserver. Since $-i\varphi(i\sigma_x)$ and $-i\varphi(i\sigma_y)$ are anticommuting self-adjoint unitaries, the dimension of the Hilbert space $\CK$ must be even whenever $\dim \CK < \infty$. Therefore, without loss of generality, we may identify $\Clb(\CK)$ with $\Clb(\CH) \otimes M_2(\Cb)$ for some Hilbert space $\CH$, and regard $\varphi \colon M_2(\Cb) \to \Clb(\CH) \otimes M_2(\Cb)$ as a unital $\Rb$-linear unitary preserver admitting the decomposition given in Theorem \ref{thm:two_conj_part}.

\begin{lem}\label{lem:conj_ipart}
    There exist a unitary $U \in \Clb(\CH) \otimes M_2(\Cb)$ and a self-adjoint unitary $H \in \Clb(\CH)$ such that  
    \begin{align*}
        U^* \varphi(i\sigma_z) U = iI \otimes \sigma_z, \quad U^* \varphi(i\sigma_y) U = iI \otimes \sigma_y, \quad U^* \varphi(i\sigma_x) U = iH \otimes \sigma_x. 
    \end{align*}
\end{lem}

\begin{proof}
    By Lemma \ref{lem:wei2} and Theorem \ref{thm:two_conj_part}, $X = -i\varphi(i\sigma_x)$, $Y = -i\varphi(i\sigma_y)$, and $Z = -i\varphi(i\sigma_z)$ are mutually anti-commuting self-adjoint unitary operators. 
    It is clear that there exists a unitary $U \in \Clb(\CH) \otimes M_2(\Cb)$ such that 
    \begin{align*}
        U^* Z U = I \otimes \sigma_z, \qquad U^* Y U = I \otimes \sigma_y.
    \end{align*}
    Since $U^* X U$ must anti-commute with both $U^* Y U$ and $U^* Z U$, it must be of the form $H \otimes \sigma_x$ for some self-adjoint unitary $H \in \Clb(\CH)_{\sa}$. 
\end{proof}

Using Lemma \ref{lem:conj_ipart}, up to a suitable unitary conjugation, we may assume that 
\begin{align*}
    \varphi(i\sigma_x) = iX, \qquad \varphi(i\sigma_y) = iY, \qquad \varphi(i\sigma_z) = iZ,
\end{align*}
where 
\begin{align*}
    X = H \otimes \sigma_x, \qquad Y = I \otimes \sigma_y, \qquad Z = I \otimes \sigma_z
\end{align*}
for some self-adjoint unitary $H \in \Clb(\CH)$. 
Since 
\begin{align*}
    \alpha(L) = -iJ \circ \varphi(iL) \quad \text{and} \quad \beta_2(L) = -\frac{1}{2} \bigl [J, \varphi(iL) \bigr ],
\end{align*}
it follows that
\begin{equation}\label{eq:alpah_beta2_con}
    \begin{aligned}
        \alpha(aI + x \sigma_x + y \sigma_y + z \sigma_z) &= aI + J \circ \bigl ( x X + yY + zZ \bigr ), \\ 
        \beta_2(aI + x \sigma_x + y \sigma_y + z \sigma_z) &= -\frac{i}{2} \bigl [J, x X + yY + zZ \bigr ].       
    \end{aligned}
\end{equation}
It is not hard to verify that $\alpha(L^2) = \alpha(L)^2 + \beta_2(L)^2$ for every $L \in M_2(\Cb)_{\sa}$.

\begin{thm}\label{thm:classification_M_2}
Suppose that $J \in \Clb(\CH) \otimes M_2(\Cb)$ and $H \in \Clb(\CH)$ are self-adjoint unitaries. 
Let 
\begin{align*}
    X = H \otimes \sigma_x, \qquad Y = I \otimes \sigma_y, \qquad Z = I \otimes \sigma_z.
\end{align*}
For each $k \in \{x, y, z\}$, let $K$ denote the corresponding operator in $\{X, Y, Z\}$, and let $U_k$ be a self-adjoint unitary on the support of $[J, K]$. 
If 
\begin{enumerate}
    \item for each $k \in \{x, y, z\}$, $U_k$ commutes with both $K$ and $J$, and
    \item for any distinct $k, l \in \{x, y, z\}$ with corresponding operators $K, L \in \{X, Y, Z\}$, we have
    \begin{align*}
        &(U_k [J, K]) \circ L + (U_l [J, L]) \circ K = 0,\\
        &(U_k [J, K]) \circ (U_l [J, L]) = [J, K] \circ [J, L],
    \end{align*}
\end{enumerate}
then $J$ and the $\Rb$-linear maps 
\begin{align*}
    \alpha(aI + x \sigma_x + y \sigma_y + z \sigma_z) &= aI + J \circ (x X + yY + zZ), \\ 
    \beta_1(aI + x\sigma_x + y\sigma_y + z\sigma_z) &= \frac{-i}{2} [J, xU_x X + y U_yY + zU_zZ], \\
    \beta_2(aI + x \sigma_x + y \sigma_y + z \sigma_z) &= \frac{-i}{2}[J, x X + yY + zZ]
\end{align*}
define an $\Rb$-linear unital unitary preserver $M_2(\Cb) \to \Clb(\CH) \otimes M_2(\Cb)$ as in Theorem \ref{thm:two_conj_part}. 
Conversely, up to unitary equivalence, every $\Rb$-linear unital unitary preserver $M_2(\Cb) \to \Clb(\CH) \otimes M_2(\Cb)$ admits the above form.
\end{thm}

\begin{proof}
    Based on the proceding discussion, we may assume that $\alpha$ and $\beta_2$ take the form given in \eqref{eq:alpah_beta2_con}. For $k \in \{x, y, z\}$, by Theorem \ref{thm:two_conj_part}, $\beta_1(\sigma_k)$ and $\beta_2(\sigma_k)$ commute and have the same square. Thus, there exists a self-adjoint unitary operator $U_k$ on the support of $\beta_2(\sigma_k)$ such that $U_k\beta_2(\sigma_k) = \beta_2(\sigma_k)U_k$ and $\beta_1(\sigma_k) = U_k \beta_2(\sigma_k)$. Since $J$ anti-commutes with both $\beta_1(\sigma_k)$ and $\beta_2(\sigma_k)$, we have   
   
    \begin{align*}
        U_k \beta_2(\sigma_k) = \beta_1(\sigma_k) = - J \beta_1(\sigma_k) J = - J U_k \beta_2(\sigma_k) J = J U_k J \beta_2(\sigma_k).
    \end{align*}
    It follows that $J U_k = U_k J$. 
    Note that 
    \begin{align*}
       0 = [U_k \beta_2(\sigma_k), \alpha(\sigma_k)] = [U_k, \alpha(\sigma_k)]\beta_2(\sigma_k).
    \end{align*}
    We have $U_k$ commutes with $\alpha(\sigma_k)$. Combining with $JU_k = U_kJ$, we conclude that  
    \begin{align*}
        U_x X = X U_x, \qquad U_y Y = Y U_y, \qquad U_z Z = Z U_z.
    \end{align*}
    Then, it is not hard to verify that  
    \begin{align*}
        \beta_1(aI + x\sigma_x + y \sigma_y + z\sigma_z) = \frac{-i}{2} [J, xU_x X + y U_yY + zU_zZ]
    \end{align*}
    for every $a, x, y, z \in \Rb$. 

    Recall that two trace $0$ self-adjoint operators in $M_2(\Cb)$ commute if and only if they are collinear. 
    Then the commutativity condition (III) in Theorem~\ref{thm:two_conj_part} amounts to the fact that $[J, xU_x X + y U_yY + zU_zZ]$ commutes with both
    \begin{align*}
        J \circ (xX + yY + zZ) \quad \text{and} \quad [J, xX + yY + zZ]. 
    \end{align*}
    Since $J$ and $[J, xU_x X + y U_yY + zU_zZ]$ anticommute, the above condition is equivalent to
    \begin{align*}
        0 =& [J, xU_x X + y U_yY + zU_zZ] \circ (xX + yY +zZ) \\
        = & xy \Bigl( (U_y [J, Y])\circ X +  (U_x [J, X]) \circ Y \Bigr) \\ 
        & + yz \Bigl( (U_z [J, Z]) \circ Y + (U_y [J, Y]) \circ Z \Bigr) \\ 
        & + zx \Bigl( (U_x [J, X]) \circ Z + (U_z [J, Z])\circ X  \Bigr),
    \end{align*}
    which yields $(U_k [J, K]) \circ L + (U_l [J, L]) \circ K = 0$ for all distinct $k, l \in \{x, y, z\}$. 
    Finally, the condition $\beta_1(x\sigma_x + y\sigma_y + z\sigma_z)^2 = \beta_2(x\sigma_x + y\sigma_y + z\sigma_z)^2$ is equivalent to 
    \begin{align*}
        (U_k [J, K]) \circ (U_l [J, L]) = [J, K] \circ [J, L].
    \end{align*}
\end{proof}

\begin{rem}
Observe that for any self-adjoint unitary $U \in \{J, X, Y, Z\}'$, setting $U_x = U_y = U_z = U$ satisfies the conditions in Theorem \ref{thm:classification_M_2}. In general, however, the operators $U_x$, $U_y$ and $U_z$ need not to be equal.
\end{rem}

\begin{cor}
Let $\varphi: M_2(\Cb) \to M_2(\Cb)$ be a unital $\Rb$-linear unitary preserver. Then $\varphi$ is either of pure Jordan-type or purely non-Jordan. If $\varphi$ is of pure Jordan-type, then $\varphi$ is unitarily equivalent to exactly one of the following maps:
\begin{itemize}
    \item $\varphi(A) = A$,
    \item $\varphi(A) = A^T$,
    \item $\varphi(A) = A^*$,
    \item $\varphi(A) = \overline{A}$.
\end{itemize}
And $\varphi$ is purely non-Jordan if and only if $\tau(\varphi(iI)) = 0$, where $\tau$ is the normalized trace on $M_2(\Cb)$. In this case, there exist $\vec{J} \in (\Rb^3)_1$ and $h,s \in\{1,-1\}$ such that $\varphi$ is unitarily equivalent to 
\begin{align*}
 aI+\vec{v} \cdot \vec{\sigma} +i\bigl(bI+\vec{u} \cdot\vec{\sigma} \bigr) \mapsto
 \bigl(a+\vec{J} \cdot D_h\vec{v} \bigr)I
 +i\Bigl( s\vec{J} \times (D_h \vec{v}) + b \vec{J} + D_h \vec{u} \Bigr) \cdot \vec \sigma,
\end{align*}
where
\begin{align*}
    \vec{\sigma} = (\sigma_x, \sigma_y, \sigma_z), \qquad D_h(v_x,v_y,v_z)=(hv_x,v_y,v_z),
\end{align*}
$a, b \in \Rb$, $\vec{u}, \vec{v} \in \Rb^3$.
\end{cor}

\begin{proof}
If $\varphi(iI) = \pm iI$, then by Corollary \ref{cor:beta_zero} and \cite[Theorem 3.3]{stormer65}, $\varphi$ is of pure Jordan type and is unitarily equivalent to one of the maps listed above.
Assume now that $\tau(\varphi(iI)) = 0$. 
By Theorem~\ref{thm:classification_M_2}, we may assume that 
\begin{align*}
    \varphi\left(i(bI + \vec{u}\cdot \vec{\sigma})\right) = i\left( b\vec{J} + D_h\vec{u} \right) \cdot \vec{\sigma}, 
\end{align*}
where $\vec{J} \in (\Rb^3)_1$. In particular, $J:= -i\varphi(iI) = \vec{J} \cdot \vec{\sigma}$. Since
\begin{align*}
 (\vec{J} \cdot\vec{\sigma})\circ(\vec{u} \cdot\vec{\sigma}) =(\vec{J}\cdot\vec{u})I,\qquad 
 [\vec{J} \cdot \vec{\sigma},\vec{u} \cdot \vec{\sigma}] =2i(\vec{J} \times\vec{u})\cdot \vec{\sigma},
\end{align*}
we have 
\begin{align*}
 \alpha(bI + \vec{u}\cdot \vec{\sigma}) = (b + \vec{J} \cdot (D_h\vec{u}))I,\qquad 
 \beta_2(bI + \vec{u}\cdot \vec{\sigma}) = (\vec{J} \times (D_h\vec{u}))\cdot\vec{\sigma}.
\end{align*}
 
Let $K \in \{h \sigma_x, \sigma_y, \sigma_z \}$. Then $[J, K]$ is either $0$ or a non-zero multiple of a traceless self-adjoint matrix that is not collinear with $J$. By condition (2) of Theorem~\ref{thm:classification_M_2}, there exists $s \in \{1,-1\}$ such that 
\begin{align*}
    \beta_1(aI+ \vec{v}\cdot \vec{\sigma}) = s (\vec{J} \times (D_h\vec{v})) \cdot \vec{\sigma}.
\end{align*}
Combining the formulas for $\alpha$, $\beta_1$, and $\beta_2$ completes the proof.
\end{proof}

\section{Binary $\Rb$-linear unital unitary preservers}

In this section, we assume that $\Clc$ is a unital C$^*$-algebra that has no irreducible representations on $\Cb^2$, or equivalently, that its enveloping von Neumann algebra $\Clc^{**}$ contains no direct summand of type $\mathrm{I}_2$, and $\varphi: \Clc \to \Clb(\CH)$ is a \textbf{purely non-Jordan} $\Rb$-linear unital unitary preserver defined by 
\begin{align}\label{equ:for_same_beta}
    \varphi(H+iK) = \alpha(H)+ i\beta(H) + iJ(\alpha(K) + i\beta(K)), 
\end{align}
where $J$ is a self-adjoint unitary in $\Clb(\CH)$, and $\alpha, \beta: \Clc_{\sa} \to \Clb(\CH)_{\sa}$ are $\Rb$-linear maps that satisfy the hypotheses of Theorem \ref{thm:two_conj_part}. We refer to these maps as binary unitary preservers. 

\begin{rem}
Let $\Cla$ be a von Neumann algebra with no direct summand of type $\mathrm{I}_2$. Then $\Cla^{**}$ also has no direct summand of type $\mathrm{I}_2$. Indeed, the kernel of every $*$-homomorphism from $\Cla$ to $M_2(\Cb)$ contains the orthogonal complement of the central projection corresponding to the abelian summand of $\Cla$, and hence the image of any such map is necessarily abelian.
\end{rem}

Let $\widehat{\varphi}: \Clc^{**} \to \Clb(\CH)$ be the normal extension of $\varphi$. For all $H, K \in (\Clc^{**})_{\sa}$, we have
\begin{align*}
    \widehat{\varphi}(H + i K) = \widehat{\alpha}(H) + i\widehat{\beta}(H) + iJ\bigl(\widehat{\alpha}(K) + i\widehat{\beta}(K)\bigr),
\end{align*}
where 
\begin{align*}
    \widehat{\alpha}(H) = \frac{1}{2} \bigl(\widehat{\varphi}(H) + \widehat{\varphi}(H)^*\bigr), \qquad 
    \widehat{\beta}(H) = \frac{1}{2i} \bigl(\widehat{\varphi}(H) - \widehat{\varphi}(H)^*\bigr) 
\end{align*}
are the normal extensions of $\alpha$ and $\beta$, respectively. Let $\Clw^*(\alpha)$ and $\Clw^*(\alpha, \beta)$ be the von Neumann algebras generated by the range of $\alpha$ and by the union of the ranges of $\alpha$ and $\beta$, respectively. Then $\Clw^*(\widehat{\alpha}) = \Clw^*(\alpha)$ and $\Clw^*(\widehat{\alpha}, \widehat{\beta}) = \Clw^*(\alpha, \beta)$. By Remark \ref{rem:ext_non_Jordan}, $\widehat{\varphi}$ is also a binary unitary preserver. In light of this observation, we hereafter restrict our attention to the case where $\Cla$ is a von Neumann algebra with no type $\mathrm{I}_2$ direct summand and $\varphi: \Cla \to \Clb(\CH)$ is a normal, $\Rb$-linear, unital binary unitary preserver. The corresponding result for $C^*$-algebras then follows by restriction. We note that the normality of $\varphi$ is required only in the proof of Theorem~\ref{thm:cls_binary}. \\

We begin our discussion by presenting a canonical example of a binary unitary preserver.

\begin{example}\label{exam:binary_exp}
Let $\Cla$ be a type II$_1$ factor, and let $(\Cla \overline{\otimes} \Cla) \rtimes \Zb_2$ be the crossed product induced by the flip $*$-automorphism of $\Cla \overline{\otimes} \Cla$. Then the self-adjoint unitary $J$ implementing the flip automorphism, together with the maps $\alpha(H) = \frac{1}{2}(H \otimes I + I \otimes H)$ and $\beta(H) = \frac{1}{2}(H \otimes I - I \otimes H)$, defines a binary unitary preserver. In particular, if $\Cla$ contains a von Neumann subalgebra isomorphic to $(\Cla \overline{\otimes} \Cla) \rtimes \Zb_2$, then composing this map with the embedding $(\Cla \overline{\otimes} \Cla) \rtimes \Zb_2 \hookrightarrow \Cla$ yields a (non-Jordan) unitary preserver on $\Cla$. Examples of such algebras include the group von Neumann algebras of infinite discrete ICC groups satisfying $G \simeq G \wr \Zb_2$; in particular, the hyperfinite II$_1$ factor belongs to this class.
\end{example}

The remainder of this section shows that binary unitary preservers essentially take the form given above. To this end, we rely on the following technical theorem, whose proof is deferred to the appendix.

\begin{thm}\label{thm:decmp_general}
Let $\Clc$ be a unital C$^*$-algebra such that $\Clc^{**}$ contains no type I$_2$ direct summand, and let $\phi_1, \phi_2: \Clc \to \Clb(\CH)$ be two unital Jordan $*$-homomorphisms such that
\begin{align*}
    [\phi_1(H), \phi_2(K)] = 0   
\end{align*} 
for all commuting operators $H, K \in \Clc_{\sa}$. Then there is a projection $P$ commuting with both ranges such that
   \begin{align*}
       \phi_1(A_1)P = \phi_2(A_1)P, \qquad [(I-P)\phi_1(A_1), (I-P)\phi_2(A_2)] = 0, \qquad \forall A_1, A_2 \in \Clc.
   \end{align*}
\end{thm}

For every projection $P$ in the set of projections $\CP(\Cla)$ in $\Cla$, we use $S_P$ to denote the support projection of $\beta(P)$, and $V_P$ to denote the self-adjoint partial isometry in the polar decomposition of $\beta(P)$, i.e., $\beta(P) = V_P |\beta(P)|$.

\begin{lem}\label{lem:bin_Jordan}
There exist unital Jordan $*$-homomorphisms $\phi_{+}, \phi_{-}: \Cla \to \Clw^*(\alpha, \beta)$ such that 
    \begin{equation}\label{equ:bin_Jordan_I}
     \begin{aligned}
        \phi_{+}(P) = \alpha(P)(I-S_P) + \frac{1}{2}(S_P + V_P), \qquad \phi_{-}(P) = \alpha(P)(I-S_P) + \frac{1}{2}(S_P - V_P)
    \end{aligned}       
    \end{equation}
    for every $P \in \CP(\Cla)$. 
    Moreover, we have the following conclusions:
    \begin{enumerate}
        \item[I.] $J\phi_{+}(A)J = \phi_{-}(A)$ for all $A \in \Cla$;
        \item[II.] the ranges of $\phi_{+}$ and $\phi_{-}$ mutually commute;  
        \item[III.] For all commuting operators $H, K \in \Cla_{\sa}$, 
        \begin{align*}
            [\alpha(H), \phi_{\pm}(K)] = [\beta(H), \phi_{\pm}(K)] = 0.
        \end{align*}
    \end{enumerate}
\end{lem}

\begin{proof}
For every $P \in \CP(\Cla)$, let 
\begin{align*}
    \rho_{\pm}(P) := \alpha(P)(I-S_P) + \frac{1}{2}(S_P \pm V_P).
\end{align*}
Let $P_1, \ldots, P_m \in \CP(\Cla)$ be a family of mutually orthogonal projections with sum $I$. We first show that $\{\rho_{+}(P_j)\}_{j=1}^{m}$ and $\{\rho_{-}(P_j)\}_{j=1}^m$ are families of mutually orthogonal projections such that  
\begin{align*}
 \sum_{j=1}^m \rho_{+}(P_j) = I = \sum_{j=1}^m \rho_{-}(P_j). 
\end{align*}

By Theorem \ref{thm:two_conj_part}, the von Neumann algebra generated by $\alpha(P_j), \beta(P_j)$, $j = 1, \ldots, m$, is abelian. By identifying the von Neumann algebra with $C(\Omega)$, where $\Omega$ is a hyperstonean compact Hausdorff space, we may treat $\alpha(P_j)$, $\beta(P_j)$, and $S_{P_j}$ as functions on $\Omega$. For each $j$, we have $\{\omega\in\Omega:S_{P_j}(\omega)=1\}=\overline{\{\omega\in\Omega:\beta(P_j)(\omega)\neq 0\}}$. Therefore, the set
\begin{align*}
 \Omega_j := \{\omega\in\Omega:S_{P_j}(\omega)=0\} \cup \{\omega\in\Omega:\beta(P_j)(\omega)\neq 0\}
\end{align*}
is dense and open in $\Omega$. 
Moreover, $\Omega_0:=\bigcap_{j=1}^m\Omega_j$ is dense and open in $\Omega$. For every $\omega\in\Omega_0$, we have $S_{P_j}(\omega)=1$ if and only if $\beta(P_j)(\omega)\neq 0$.

Fix $\omega \in \Omega_0$. Since $\alpha(I) = I$, $\beta(I) = 0$, and $\alpha(P_j) \alpha(P_k) + \beta(P_j)\beta(P_k) = \delta_{j,k}\alpha(P_j)$, the vectors $(\alpha(P_j)(\omega),\beta(P_j)(\omega)) \in\Rb^2$ for $j=1,\ldots,m$ are mutually orthogonal, and
\begin{align*}
 \sum_j\left(\alpha(P_j)(\omega),\beta(P_j)(\omega)\right) =(1,0), \qquad
 \alpha(P_j)(\omega)^2+\beta(P_j)(\omega)^2 = \alpha(P_j)(\omega).
\end{align*}
Since $\Rb^2$ has dimension $2$, there can be at most two mutually orthogonal nonzero vectors.
If $(\alpha(P_j)(\omega),\beta(P_j)(\omega)) = (1, 0)$ is the only nonzero vector, then 
\begin{align*}
    \rho_{\pm}(P_k)(\omega) = 
    \begin{cases}
        1 & k = j\\
        0 & k \neq j     
    \end{cases}.
\end{align*}
Assume now that the $j$-th and $k$-th vectors are nonzero. Without loss of generality, we may assume that 
\begin{align*}
 \left(\alpha(P_j)(\omega),\beta(P_j)(\omega)\right) = (a, \sqrt{a(1-a)}), \qquad 
 \left(\alpha(P_k)(\omega),\beta(P_k)(\omega)\right) = (1-a,- \sqrt{a(1-a)}),
\end{align*}
where $0<a<1$. Note that $\beta(P)=V_P|\beta(P)|$. We have 
\begin{align*}
     \rho_{+}(P_l)(\omega) = 
    \begin{cases}
        1 & l = j\\
        0 & l \neq j     
    \end{cases}, \qquad
     \rho_{-}(P_l)(\omega) = 
    \begin{cases}
        1 & l = k\\
        0 & l \neq k    
    \end{cases}.
\end{align*}
Therefore, there is a unique index $j$ such that $\rho_{+}(P_j)(\omega) = 1$ (resp., $\rho_{-}(P_j)(\omega) = 1$) and $\rho_{+}(P_k)(\omega) = 0$ (resp., $\rho_{-}(P_k)(\omega) = 0$) for all $k \neq j$. Since $\Omega_0$ is dense in $\Omega$, we conclude that $\rho_{\pm}(P_1), \ldots, \rho_{\pm}(P_m)$ are mutually orthogonal projections with sum $I$. 
In particular, 
\begin{align*}
    \rho_{\pm}(P+Q) = \rho_{\pm}(P) + \rho_{\pm}(Q) = I -\rho_{\pm}(I-P - Q)
\end{align*}
for every pair of mutually orthogonal projections $P, Q$. By Theorem A in \cite{bunce92}, $\rho_{\pm}$ uniquely extends to a bounded unital $\Cb$-linear map $\phi_{\pm}$ from $\Cla$ to $\Clw^*(\alpha, \beta)$. 
Note that $\phi_{\pm}$ maps mutually orthogonal projections to mutually orthogonal projections. It is routine to check that $\phi_{+}$ and $\phi_{-}$ are Jordan $*$-homomorphisms (see, for example, Theorem A.4 in \cite{M2007}).

Recall that $J\beta(P)J = -\beta(P)$, which implies $JS_PJ = S_P$ and $JV_PJ = - V_P$ for every projection $P \in \Cla$. By approximating self-adjoint operators with linear combinations of mutually orthogonal projections and applying Equation \eqref{equ:bin_Jordan_I}, it is routine to check that $J\phi_{+}(A)J = \phi_{-}(A)$ for all $A \in \Cla$, and
\begin{align*}
  [\phi_{+}(H), \phi_{-}(K)] = [\alpha(H), \phi_{\pm}(K)] = [\beta(H), \phi_{\pm}(K)] = 0
\end{align*}
whenever $H, K \in \Cla_{\sa}$ commute. 

By Theorem \ref{thm:decmp_general}, there exists a projection $F \in \Clb(\CH)$ commuting with the ranges of $\phi_{+}$ and $\phi_{-}$ such that 
\begin{align*}
   \phi_{+}(A)F = \phi_{-}(A)F, \qquad [\phi_{+}(A)(I-F), \phi_{-}(B)(I-F)] = 0, \qquad \forall A, B \in \Cla. 
\end{align*}
It remains to show that $F =0$. For every $P \in \CP(\Cla)$, we have $FV_P = F(\phi_{+}(P) - \phi_{-}(P)) = 0$. Therefore
\begin{align*}
    F\beta(P) = \beta(P)F = 0, \qquad F\alpha(P) = F\phi_{+}(P) = \phi_{+}(P)F = \alpha(P)F, \qquad \forall P \in \CP(\Cla).
\end{align*}
By Proposition \ref{prop:Jordan_center_split}, the central support of $S_{\varphi} = \bigvee_{P \in \CP(\Cla)} S_P$ in $\Clw^*(\alpha)$ is $I$. 
Since $F \alpha(H) S_{\varphi} \xi =  \alpha(H) F S_{\varphi} \xi = 0$ for all $\xi \in \CH$, we obtain $F =0$.  
\end{proof}

\begin{lem}\label{lem:Jordan-completion}
Let
\begin{align}\label{equ:Jordan-completion_0}
    \gamma(A):= \phi_{+}(A) + \phi_{-}(A) - \alpha(A), \qquad A \in \Cla.
\end{align}
Then the map
\begin{align}\label{equ:Jordan-completion_I}
    \Pi(A) :=
    \begin{pmatrix}
        \alpha(A) & \beta(A)\\
        \beta(A) & \gamma(A)
    \end{pmatrix},
    \qquad A\in\Cla
\end{align}
is a unital Jordan $*$-homomorphism from $\Cla$ into the commutant of $\diag(J, -J)$ in $\Clb(\CH) \otimes M_2(\Cb)$, where $\alpha,\beta$ denote the $\Cb$-linear extensions of the original maps.
\end{lem}

\begin{proof}
For any projection $P\in\CP(\Cla)$, we have $\gamma(P)=\alpha(P)(I-S_P)+(I-\alpha(P))S_P$. Using the commutativity of $\alpha(P)$ and $\beta(P)$ along with the relation $\alpha(P) = \alpha(P)^2 + \beta(P)^2$, it is straightforward to check that $\Pi(P)^2 = \Pi(P)$, whence $\Pi(P)$ is a projection. Then $\Pi$ is a unital Jordan $*$-homomorphism by Theorem A.4 in \cite{M2007}. By Lemma \ref{lem:bin_Jordan}, $J \gamma(A) J = \gamma(A)$. 
Combined with the identity $J\beta(A)J = -\beta(A)$ for all $A \in \Cla$, we obtain that $\Pi(A)$ commutes with $\diag(J, -J)$.
\end{proof}

The unital Jordan $*$-homomorphism $\Pi \colon \Cla \to \Clb(\CH) \otimes M_2(\Cb)$ defined by \eqref{equ:Jordan-completion_I} will be referred to as the canonical Jordan $*$-homomorphism associated with the $\Rb$-linear unital unitary preserver $\varphi$.

\begin{lem}\label{lem:diag_unique}
If a positive contraction $T$ commuting with the ranges of $\alpha$ and $\beta$ satisfies
\begin{align*}
   \Pi(A) =  
   \begin{pmatrix}
      \sqrt{T} & \sqrt{I-T}\\
      \sqrt{I-T} & -\sqrt{T}
   \end{pmatrix} 
   \begin{pmatrix}
    \phi_{1}(A) & 0\\
    0 & \phi_{2}(A)
    \end{pmatrix}
   \begin{pmatrix}
      \sqrt{T} & \sqrt{I-T}\\
      \sqrt{I-T} & -\sqrt{T}
   \end{pmatrix} 
\end{align*}
for some Jordan $*$-homomorphisms $\phi_{1}$ and $\phi_{2}$, then $T$ is necessarily unique and belongs to the center of $\Clw^*(\alpha, \beta)$. Moreover, $\phi_1 = \phi_{+}$, $\phi_2 = \phi_{-}$, $JTJ = I-T$, and $\supp T(I-T) = I$.
\end{lem}

\begin{proof} 
Assume that 
\begin{align}\label{equ:diag_unique_I}
   \begin{pmatrix}
      \sqrt{T} & \sqrt{I-T}\\
      \sqrt{I-T} & -\sqrt{T}
   \end{pmatrix} 
   \begin{pmatrix}
    \alpha(A) & \beta(A)\\
    \beta(A) & \gamma(A)
    \end{pmatrix}
   \begin{pmatrix}
      \sqrt{T} & \sqrt{I-T}\\
      \sqrt{I-T} & -\sqrt{T}
   \end{pmatrix} 
   =
   \begin{pmatrix}
    \phi_{1}(A) & 0\\
    0 & \phi_{2}(A)
    \end{pmatrix}.
\end{align}
First note that $T$ also commutes with the range of $\gamma$ by Equation \eqref{equ:bin_Jordan_I}. Let $P \in \CP(\Cla)$. Note that $\alpha(P) - \gamma(P) = S_P(2\alpha(P)-I)$. 
Then 
\begin{align}\label{equ:diag_unique_II}
  (\alpha(P) - \gamma(P))^2 + 4\beta(P)^2 = S_P(2\alpha(P)-I)^2 + 4\beta(P)^2 = S_P.
\end{align}
By comparing the $(1, 2)$-entries on both sides of Equation \eqref{equ:diag_unique_I}, we obtain 
\begin{align}\label{equ:diag_unique_III}
 \sqrt{T(I-T)}(\alpha(P)-\gamma(P)) = (2T-I)\beta(P).
\end{align}
Squaring both sides of the Equation \eqref{equ:diag_unique_III} and combining with Equation \eqref{equ:diag_unique_II}, we have 
\begin{align*}
  T(I-T) S_P = T(I-T)\left [ (\alpha(P)-\gamma(P))^2 + 4\beta(P)^2 \right ] = \beta(P)^2 .
\end{align*}
It follows that  
\begin{align}\label{eq:diag_unique_V}
  \sqrt{T(I-T)}S_P = |\beta(P)|.  
\end{align}
By Equations \eqref{equ:diag_unique_III} and \eqref{eq:diag_unique_V}, we have 
\begin{align*}
   (2T - I)\beta(P) = \sqrt{T(I-T)} S_P(2\alpha(P)-I) = (2\alpha(P)-I)|\beta(P)|.
\end{align*}
Hence,
\begin{align}\label{equ:diag_unique_IV}
    (2T - I)S_P = (2\alpha(P)-I)V_P, \qquad \forall P \in \CP(\Cla).
\end{align}
Let $M$ be an operator commuting with the ranges of $\alpha$ and $\beta$. By Equation \eqref{equ:diag_unique_IV}, 
\begin{align*}
    S_PTM &= \frac{1}{2}[(2\alpha(P)-I)V_P + S_P]M = S_PMT, \quad P \in \CP(\Cla).
\end{align*}
Since the central support of $\bigvee_{P \in \CP(\Cla)} S_P$ in $\Clw^*(\alpha)$ is $I$, we have $TM = MT$. 
Therefore, $T$ lies in the center of $\Clw^*(\alpha, \beta)$. Similarly, Equations \eqref{eq:diag_unique_V} and \eqref{equ:diag_unique_IV}, together with the fact that the central support of $\bigvee_{P \in \CP(\Cla)} S_P$ is $I$, ensure the uniqueness of $T$ and that $\supp T(I-T) = I$.

Comparing the $(1,1)$-entry and the $(2,2)$-entry of Equation \eqref{equ:diag_unique_I}, we have
\begin{align*}
    \phi_1(P) &= T\alpha(P) + (I-T)\gamma(P) + 2\sqrt{T(I-T)}\beta(P),\\
    \phi_2(P) &= (I-T)\alpha(P) + T\gamma(P) - 2\sqrt{T(I-T)}\beta(P).
\end{align*}
By Equations \eqref{eq:diag_unique_V} and \eqref{equ:diag_unique_IV},  
\begin{align*}
    \phi_1(P) + \phi_2(P) &= \alpha(P) + \gamma(P) = \phi_{+}(P) + \phi_{-}(P), \\
    \phi_1(P) - \phi_2(P) &= (2T-I)(\alpha(P) - \gamma(P)) + 4\sqrt{T(I-T)}\beta(P) \\
    &= (2T-I)^2 V_P + 4T(I-T)V_P = V_P = \phi_{+}(P) - \phi_{-}(P).  
\end{align*}
This implies that $\phi_1 = \phi_{+}$ and $\phi_2 = \phi_{-}$.

Recall from Lemma \ref{lem:Jordan-completion} that $\Pi(A) = \diag(J, -J)\Pi(A) \diag(J, -J)$. 
This yields 
\begin{align*}
    \Pi(A) &=  
   \begin{pmatrix}
      \sqrt{JTJ} & -\sqrt{J(I-T)J}\\
      -\sqrt{J(I-T)J} & -\sqrt{JTJ}
   \end{pmatrix} 
   \begin{pmatrix}
    \phi_{-}(A) & 0\\
    0 & \phi_{+}(A)
    \end{pmatrix}
   \begin{pmatrix}
      \sqrt{JTJ} & -\sqrt{J(I-T)J}\\
      -\sqrt{J(I-T)J} & -\sqrt{JTJ}
   \end{pmatrix}\\ 
   &=    \begin{pmatrix}
      \sqrt{J(I-T)J} & \sqrt{JTJ}\\
      \sqrt{JTJ} & -\sqrt{J(I-T)J}
   \end{pmatrix} 
   \begin{pmatrix}
    \phi_{+}(A) & 0\\
    0 & \phi_{-}(A)
    \end{pmatrix}
   \begin{pmatrix}
      \sqrt{J(I-T)J} & \sqrt{JTJ}\\
      \sqrt{JTJ} & -\sqrt{J(I-T)J}
   \end{pmatrix}. 
\end{align*}
The uniqueness of $T$ then implies that $JTJ = I-T$.
\end{proof}

In the remainder of this section, we establish the existence of a positive contraction $T$ satisfying the properties given in Lemma \ref{lem:diag_unique}. The following remark clarifies why the type I$_2$ direct summand must be excluded from our assumptions.

\begin{rem}\label{rem:I2_counter_exp}
By Theorem \ref{thm:classification_M_2}, up to unitary equivalence, every $\Rb$-linear unital unitary preserver
\begin{align*}
\varphi:M_2(\Cb)\longrightarrow \Clb(\CH)\otimes M_2(\Cb)
\end{align*}
for which the decomposition in Theorem \ref{thm:two_conj_part} satisfies $\beta_1=\beta_2=\beta$ admits the following form:
\begin{align*}
\varphi(aI+x\sigma_x+y\sigma_y+z\sigma_z) &= aI+J(xX+yY+zZ),\\
\varphi\left(i(aI+x\sigma_x+y\sigma_y+z\sigma_z)\right) &=iaJ+i(xX+yY+zZ),
\end{align*}
where $J\in\Clb(\CH)\otimes M_2(\Cb)$ and $H\in\Clb(\CH)$ are self-adjoint unitaries, and
\begin{align*}
 X=H\otimes\sigma_x,\qquad Y=I\otimes\sigma_y,\qquad Z=I\otimes\sigma_z.   
\end{align*}

Setting $\CH = \Cb$, $H= 1$, and $J = \sigma_z$, we have
\begin{align*}
    \alpha(aI+x\sigma_x+y\sigma_y+z\sigma_z) = (a + z)I,\qquad \beta(aI+x\sigma_x+y\sigma_y+z\sigma_z) = x \sigma_y - y \sigma_x.
\end{align*}
In particular, $\Clw^*(\alpha, \beta) = M_2(\Cb)$. 
If there exists a $T \in M_2(\Cb)$ satisfying the conditions in Lemma \ref{lem:diag_unique}, then $T$ must equal $\frac{1}{2}I$. 
Let $\phi_1$ and $\phi_2$ be the Jordan $*$-homomorphisms such that
\begin{align*}
    \alpha(A) = \frac{1}{2}(\phi_1(A) + \phi_2(A)), \qquad \beta(A) = \frac{1}{2}(\phi_1(A) - \phi_2(A)).
\end{align*}
Then,
\begin{align*}
    \phi_1(\sigma_z) = \alpha(\sigma_z) + \beta(\sigma_z) = I,\qquad \phi_1(\sigma_x)= \alpha(\sigma_x) + \beta(\sigma_x) = \sigma_y. 
\end{align*}
However, 
\begin{align*}
    \sigma_y = \phi_1(\sigma_z) \circ \phi_1(\sigma_x) \neq \phi_1(\sigma_z \circ \sigma_x) = 0.
\end{align*}
This contradicts the fact that $\phi_1$ is a Jordan $*$-homomorphism. 
Therefore, no such $T$ exists.
\end{rem}

To prove the existence of a positive contraction $T$ satisfying the properties given in Lemma \ref{lem:diag_unique}, we divide the argument into three cases. We first consider the case where $\Cla$ is abelian.

\begin{lem}\label{lem:I1_case}
Assume that $\Cla$ is abelian. Then there exists a positive contraction $T \in \Clb(\CH)$ commuting with the ranges of $\alpha$ and $\beta$ such that
    \begin{align*}
        \Pi(A) = 
        \begin{pmatrix}
            \sqrt{T} & \sqrt{I-T}\\
            \sqrt{I-T} & -\sqrt{T}    
        \end{pmatrix}
        \begin{pmatrix}
            \phi_{+}(A) & \\
             & \phi_{-}(A)    
        \end{pmatrix}
        \begin{pmatrix}
            \sqrt{T} & \sqrt{I-T}\\
            \sqrt{I-T} & -\sqrt{T}    
        \end{pmatrix}, \quad \forall A \in \Cla,
    \end{align*}
    where $\Pi$ is the canonical Jordan $*$-homomorphism associated with $\varphi$.
\end{lem}

\begin{proof}
    Since $\alpha(P)=\alpha(P)^2+\beta(P)^2$ for $P \in \CP(\Cla)$, we have $-I \leq 2\alpha(P) - I \leq I$. Therefore, 
    \begin{align*}
    T_P := \frac{1}{2} \left(S_P + (2\alpha(P)-I) V_P\right)
    \end{align*}
    is a positive contraction on $S_P \CH$. 
    We first prove that
    \begin{align*}
    T_{E} S_{F} =T_FS_E, \qquad \forall E,F \in \CP(\Cla). 
    \end{align*}

    Let $E = P_1 + P_2$ and $F= P_1 + P_3$, where $P_1, P_2, P_3$ are mutually orthogonal projections in $\Cla$. Set $P_4 = I-P_1 - P_2- P_3$. By identifying the abelian von Neumann algebra $\Clw^*(\alpha, \beta)$ with $C(\Omega)$ for a hyperstonean compact Hausdorff space $\Omega$, we treat operators in $\Clw^*(\alpha, \beta)$ as continuous functions on $\Omega$. 
    As in the proof of Lemma \ref{lem:bin_Jordan}, we choose a dense subset $\Omega_0$ of $\Omega$ such that for every $\omega \in \Omega_0$, $S_{Q} (\omega) = 1$ if and only if $\beta(Q)(\omega) \neq 0$, for every $Q \in \{P_1, \ldots, P_4, E, F\}$. Recall that $\{(\alpha(P_j)(\omega), \beta(P_j)(\omega))\}_{j=1}^{4}$ are orthogonal vectors in $\Rb^2$ and
    \begin{align}\label{equ:lem:I1_case_I}
        \sum_{j=1}^4 (\alpha(P_j)(\omega), \beta(P_j)(\omega)) = (1, 0).
    \end{align} 
    In particular, at most two of the four vectors are nonzero.

    Since $\Omega_0$ is dense in $\Omega$ and 
    \begin{align*}
        T_E S_F = T_E S_E S_F, \qquad T_F S_E = T_F S_E S_F,
    \end{align*}
    it suffices to verify that $T_E(\omega)S_F(\omega) = T_F(\omega)S_E(\omega)$ for all $\omega \in \Omega_0$ such that $S_E(\omega)S_F(\omega) \neq 0$. 

    Fix $\omega \in \Omega_0$ satisfying $S_E(\omega)S_F(\omega) \neq 0$, so that
    \begin{align*}
        \beta(P_1)(\omega) + \beta(P_2)(\omega) \neq 0, \qquad \beta(P_1)(\omega) + \beta(P_3)(\omega) \neq 0.
    \end{align*}
    If $\beta(P_1)(\omega) \neq 0$, then $\alpha(P_2)(\omega) = \alpha(P_3)(\omega) = 0$, as otherwise Equation \eqref{equ:lem:I1_case_I} would imply $\beta(P_1)(\omega) + \beta(P_2)(\omega) = 0$ or $\beta(P_1)(\omega) + \beta(P_3)(\omega) = 0$. Therefore, in this case, $T_E(\omega)S_F(\omega) = T_F(\omega)S_E(\omega)$, since 
    \begin{align*}
        (\alpha(E)(\omega), \beta(E)(\omega)) = (\alpha(P_1)(\omega), \beta(P_1)(\omega)) = (\alpha(F)(\omega), \beta(F)(\omega)).
    \end{align*} 
    It remains to consider the case where $\beta(P_1)(\omega) = 0$, $\beta(P_2)(\omega) \neq 0$, and $\beta(P_3)(\omega) \neq 0$. We may assume that
    \begin{align*}
    (\alpha(P_2)(\omega) , \beta(P_2)(\omega) ) = (a, b), \quad (\alpha(P_3)(\omega) , \beta(P_3)(\omega) ) = (1-a, -b), 
    \end{align*} 
    where $a \in (0, 1)$ and $b$ satisfies $b^2 = a(1-a)$. Then it is straightforward to verify that 
     \begin{align*}
        T_E(\omega)S_{F}(\omega) = \frac{1}{2}(1 + (2a-1)\sgn(b))  = T_F(\omega)S_E(\omega). 
    \end{align*}  
    
    Since $\varphi$ is purely non-Jordan, we have
    \begin{align*}
      S_{\varphi} = \bigvee_{P \in \CP(\Cla)} S_P = I \mbox{ (see Proposition \ref{prop:Jordan_center_split})}.
    \end{align*} 
    Then the family $\{T_P\}_{P \in \CP(\Cla)}$ defines a positive contraction $T \in \Clw^*(\alpha, \beta)$ such that $TS_P = T_P$ and hence 
    \begin{align}\label{equ:abelian_glue_I}
        (2T-I) V_P = (2\alpha(P)-I)S_P, \quad \forall P \in \CP(\Cla).
    \end{align} 
    Squaring Equation \eqref{equ:abelian_glue_I}, we get $T(I-T)S_P = \beta(P)^2$.
    Since $\alpha(P) - \gamma(P) = (2\alpha(P)-I)S_P$, we have
    \begin{align*}
        \sqrt{T(I-T)}(\alpha(P)- \gamma(P)) = \sqrt{T(I-T)}S_P(2\alpha(P)-I) = (2\alpha(P)-I)|\beta(P)| = (2T-I) \beta(P). 
    \end{align*}
    This implies that the off-diagonal entries of  
    \begin{align*}
    \begin{pmatrix}
      \sqrt{T} & \sqrt{I-T}\\
      \sqrt{I-T} & -\sqrt{T}
   \end{pmatrix} 
   \begin{pmatrix}
    \alpha(A) & \beta(A)\\
    \beta(A) & \gamma(A)
    \end{pmatrix}
   \begin{pmatrix}
      \sqrt{T} & \sqrt{I-T}\\
      \sqrt{I-T} & -\sqrt{T}
   \end{pmatrix}        
    \end{align*}
    are $0$. 
    The conclusion now follows from Lemma \ref{lem:diag_unique}.
\end{proof}

We now turn to the case where $\Cla$ is not abelian. 

\begin{lem}\label{lem:pm_deproj_in_center}
With the notation of Lemma \ref{lem:Jordan-completion}, let $Q$ be a central projection in the von Neumann algebra generated by the image of $\phi_{\epsilon}$ such that $Q\phi_{\epsilon}$ is a $*$-homomorphism and $(I-Q)\phi_{\epsilon}$ is a $*$-anti-homomorphism, where $\epsilon\in\{+,-\}$.
Then $Q$ is in the center of $\Clw^*(\alpha, \beta)$.    
\end{lem}

\begin{proof}
Since 
\begin{align*}
    \phi_{\pm}(P) = \alpha(P)(I-S_P) + \frac{1}{2}(S_P \pm V_P), \qquad \forall P \in \CP(\Cla),
\end{align*}
the von Neumann algebra generated by the image of $\phi_{\epsilon}$ is contained in $\Clw^*(\alpha, \beta)$ for $\epsilon \in \{+, -\}$. 
Then it suffices to show that $Q$ commutes with the ranges of $\alpha$ and $\beta$. 
By Lemma \ref{lem:bin_Jordan}, 
\begin{align*}
    \begin{pmatrix}
        \phi_{\epsilon}(H) & 0\\
        0 & \phi_{\epsilon}(H)
    \end{pmatrix} \mbox{ and }
    \begin{pmatrix}
        \alpha(K) & \beta(K)\\
        \beta(K) & \gamma(K)
    \end{pmatrix}
\end{align*}
commute for all commuting operators $H, K \in \Cla_{\sa}$. By Theorem \ref{thm:decmp_general}, $Q$ commutes with the ranges of $\alpha$ and $\beta$. 
\end{proof}

We will use the following well-known fact in the upcoming proofs. A brief sketch is provided for the reader's convenience.

\begin{lem}\label{lem:fl_fact}
Let $\Clm$ be a von Neumann algebra, $\rho: \Clm \to \Clb(\CH)$ be a unital Jordan $*$-homomorphism, and $\theta: \Clm \to \Clb(\CH)$ a norm-continuous $\Cb$-linear map. 
If 
\begin{align*}
    \theta(P) = \rho(P) \theta(P) \rho(P), \qquad \forall P \in \CP(\Clm),
\end{align*}
then $\theta(A) = \theta(I) \rho(A) = \rho(A) \theta(I)$ for all $A \in \Clm$. 
In particular, if $\theta$ is a positive map less than $\rho$, then $\theta(A) = \theta(I) \rho(A)$. 
\end{lem}

\begin{proof}
Since $\theta(P) = \rho(P) \theta(P) \rho(P)$, we have 
\begin{equation*}
    \theta(I)\rho(P) = (\theta(P) + \theta(I-P)) \rho(P) = \theta(P) =  \rho(P)(\theta(P) + \theta(I-P)) = \rho(P)\theta(I).
\end{equation*} 
By approximating operators in $\Clm$ by linear combinations of projections, we conclude that $\theta(A) = \theta(I)\rho(A) = \rho(A) \theta(I)$ for every $A \in \Clm$.
 
If $0 \leq \theta \leq \rho$, then $0 \leq \theta(P) \leq \rho(P)$ for every $P \in \CP(\Clm)$. 
This implies that $\theta(P) = \rho(P) \theta(P) \rho(P)$. The conclusion then follows from the first part.
\end{proof}

\begin{lem}\label{lem:case_CM}
    Assume that $\Cla = \Cln \otimes M_n(\Cb)$, where $\Cln$ is a von Neumann algebra and $n \geq 3$. Then there exists a positive contraction $T \in \Clb(\CH)$ commuting with the ranges of $\alpha$ and $\beta$ such that
    \begin{align}\label{equ:case_CM_0}
        \Pi(A) = 
        \begin{pmatrix}
            \sqrt{T} & \sqrt{I-T}\\
            \sqrt{I-T} & -\sqrt{T}    
        \end{pmatrix}
        \begin{pmatrix}
            \phi_{+}(A) & \\
             & \phi_{-}(A)    
        \end{pmatrix}
        \begin{pmatrix}
            \sqrt{T} & \sqrt{I-T}\\
            \sqrt{I-T} & -\sqrt{T}    
        \end{pmatrix}, \quad \forall A \in \Cla,
    \end{align}
    where $\Pi$ is the canonical Jordan $*$-homomorphism associated with $\varphi$.
\end{lem}

\begin{proof}
    Let $\phi_{+}, \phi_{-}: \Cla \to \Clb(\CH)$ be the unital Jordan $*$-homomorphisms defined by \eqref{equ:bin_Jordan_I}. By Lemma \ref{lem:pm_deproj_in_center}, and upon replacing $\Cla$ with $\Cla^{\op}$ if necessary, it suffices to prove the lemma under the assumption that $\phi_{+}$ is a $*$-homomorphism and $\phi_{-}$ is a Jordan $*$-homomorphism.

    Without loss of generality, we may assume that $\CH = \CK \otimes \Cb^n$, 
    \begin{equation*}
        \phi_{+}(N \otimes M) = \sigma(N) \otimes M, \qquad \forall N \in  \Cln, M \in M_n(\Cb), 
    \end{equation*}
    and $\phi_{-}(A) = \rho(A) \otimes I_n$ for all $A \in \Cla$, where $\sigma: \Cln \to \Clb(\CK)$ is a $*$-homomorphism, and $\rho: \Cla \to \Clb(\CK)$ is a Jordan $*$-homomorphism. 
    Note that the ranges of $\sigma$ and $\rho$ mutually commute. 
    In the following, we identify $\Clb(\CK \otimes \Cb^n) \otimes M_2(\Cb)$ with $(\Clb(\CK) \otimes M_2(\Cb)) \otimes M_n(\Cb)$. 
    Under this identification,
    \begin{align*}
    \begin{pmatrix}
        \phi_{+}(N \otimes M) & 0 \\
        0 & \phi_{+}(N \otimes M) 
    \end{pmatrix} 
    \quad \mbox{and} \quad 
    \begin{pmatrix}
        \phi_{-}(A) & 0 \\
        0 & \phi_{-}(A) 
    \end{pmatrix}
    \end{align*} 
    are identified with $\tilde{\sigma}(N) \otimes M$ and $\tilde{\rho}(A) \otimes I_n$, respectively, where  
    \begin{align*}
        \tilde{\sigma}(N) = \diag(\sigma(N), \sigma(N)) \quad \mbox{ and } \quad \tilde{\rho}(A) = \diag(\rho(A), \rho(A)).
    \end{align*}
 
    By Lemma \ref{lem:bin_Jordan}, 
    \begin{align*}
        [\Pi(H), \diag(\phi_{+}(K), \phi_{+}(K))] = 0, 
    \end{align*}  
    for all mutually commuting operators $H, K \in \Cla_{\sa}$, where $\Pi$ is the canonical Jordan $*$-homomorphism associated with $\varphi$. 
    By Theorem \ref{thm:decmp_general}, there exists a projection $F_0 \in \Clb(\CK) \otimes M_2(\Cb)$ such that the projection $F := F_0 \otimes I_n$ commutes with both ranges of $\Pi$ and $\phi_{+} \otimes I_2$, 
    \begin{align*}
        F\Pi(N \otimes M) = F_0\tilde{\sigma}(N) \otimes M, \qquad N \in \Cln, M \in M_n(\Cb),
    \end{align*}
    and 
    \begin{align*}
        [(I- F)\Pi(N_1 \otimes M_1), (I-F_0)\tilde{\sigma}(N_2) \otimes M_2] = 0  
    \end{align*}
    for $N_1, N_2 \in \Cln$, and $M_1, M_2 \in M_n(\Cb)$. 
    In particular, there exists a unital Jordan $*$-homomorphism $\Psi: \Cla \to (I-F_0)(\Clb(\CK) \otimes M_2(\Cb))(I-F_0)$ with the decomposition 
    \begin{align*}
        \Psi(A) = 
        \begin{pmatrix}
            \Psi_{11}(A) & \Psi_{12}(A) \\
            \Psi_{21}(A) & \Psi_{22}(A)
        \end{pmatrix}, 
    \end{align*}
    such that 
    \begin{align}\label{equ:case_CM_I}
        \Pi(N \otimes M) =  F_0\tilde{\sigma}(N) \otimes M + \Psi(N \otimes M) \otimes I_n, \qquad N \in \Cln, M \in M_n(\Cb).
    \end{align}
    By Equation \eqref{equ:Jordan-completion_0}, the unnormalized partial trace of $\Pi(N \otimes M)$ over the tensor factor $M_2(\Cb)$ equals $\phi_{+}(N \otimes M) + \phi_{-}(N \otimes M)$. 
    Taking the partial trace over the tensor factor $M_2(\Cb)$ of Equation \eqref{equ:case_CM_I}, we obtain
    \begin{equation}\label{equ:case_CM_2}
       \begin{aligned}
        \bigl ( (\Tr_2(F_0) -I) \sigma(N) \bigr ) \otimes M 
        = \bigl (\rho(N \otimes M) -  (\Psi_{11}(N \otimes M) + \Psi_{22}(N \otimes M)) \bigr ) \otimes I_n,
       \end{aligned}   
    \end{equation}
    where $\Tr_2: \Clb(\CK) \otimes M_2(\Cb) \to \Clb(\CK)$ is the unnormalized partial trace over the second tensor factor. By taking $M$ to be a minimal projection, \eqref{equ:case_CM_2} implies that
    \begin{align}\label{equ:case_CM_5}
        \Psi_{11}(A) + \Psi_{22}(A) = \rho(A), \qquad  \forall A \in \Cla, 
    \end{align}  
    and 
    \begin{align*}
        \Tr_2(F_0) \sigma(N) = \sigma(N), \qquad \forall N \in \Cln. 
    \end{align*}
    In particular, $\Tr_2(F_0) = I_{\CK}$. For every projection $P \in \Cla$, the Equation \eqref{equ:case_CM_5} implies 
    \begin{align*}
       \Psi(P) \leq \diag(\rho(P), \rho(P)).
    \end{align*}
    By Lemma \ref{lem:fl_fact}, we have 
    \begin{align*}
       \Psi(A) = (I-F_0) \diag(\rho(A), \rho(A)) = \diag(\rho(A), \rho(A))(I-F_0), \qquad \forall A \in \Cla. 
    \end{align*} 
    Consequently,
    \begin{align}\label{equ:case_CM_4}
       \Pi(A) = F \diag(\phi_{+}(A), \phi_{+}(A)) + (I-F)\diag(\phi_{-}(A), \phi_{-}(A)) 
    \end{align}
    for every $A \in \Cla$. 
 
    Since $\Tr_2(F_0) = I_{\CK}$, we can decompose $F$ as
    \begin{align*}
        F = 
        \begin{pmatrix} 
            T & Y \\ 
            Y^* & I-T 
        \end{pmatrix},
    \end{align*}
    where $T$ is a positive contraction, and  $Y \in \Clb(\CK \otimes \Cb^n)$. Since $F^2 = F$, we have $YY^* = T(I-T)$ and $TY = YT$. Comparing the $(1,1)$ and $(1,2)$ entries in Equation \eqref{equ:case_CM_4} gives
    \begin{align}\label{equ:case_CM_3}
         \alpha(A) = T\phi_{+}(A) + (I-T)\phi_{-}(A), \qquad \beta(A) &= Y(\phi_{+}(A) - \phi_{-}(A))
    \end{align} 
    for every $A \in \Cla$. Combined with the fact that $T$ and $Y$ commute with the ranges of $\phi_{+}$ and $\phi_{-}$, the equations in \eqref{equ:case_CM_3} imply that $T$ and $Y$ belong to the commutant of $\Clw^*(\alpha, \beta)$. Using the second equation in \eqref{equ:case_CM_3}, we obtain
    \begin{align*}
        Y S_P = (YV_P)V_P = \beta(P)V_P = |\beta(P)| \geq 0.
    \end{align*}
    In particular, this implies that $(Y - Y^*) S_P = 0$. Because the central support of 
    \begin{equation*}
        \bigvee_{P \in \CP(\Cla)} S_P = \bigvee_{P \in \CP(\Cla)} (\phi_{+}(P) - \phi_{-}(P))^2
    \end{equation*}
    in $\Clw^*(\alpha)$ is the identity $I$, it follows that $Y = Y^*$. Let $E_{n}$ be the spectral projection of $Y$ corresponding to the interval $(-\infty, -\frac{1}{n})$. We have
    \begin{equation*}
        0 \leq YE_{n} S_P \leq -\frac{1}{n} E_nS_P \leq 0,
    \end{equation*}
    which implies that $E_nS_P = 0$. 
    Since the central support of $\bigvee_{P \in \CP(\Cla)} S_P$ is $I$, it must be that $E_n = 0$ for all $n$. 
    Consequently, $Y \geq 0$, which yields $Y = \sqrt{T(I-T)}$. 
    Note that
    \begin{align*}
        F = \begin{pmatrix} 
            \sqrt{T} & \sqrt{I-T} \\
            \sqrt{I-T} & -\sqrt{T} 
        \end{pmatrix}
        \begin{pmatrix} 
            I &  \\
             & 0 
        \end{pmatrix}
        \begin{pmatrix} 
            \sqrt{T} & \sqrt{I-T} \\
            \sqrt{I-T} & -\sqrt{T} 
        \end{pmatrix}.
    \end{align*}
    Then the Equation \eqref{equ:case_CM_0} is amount to Equation \eqref{equ:case_CM_4}.
\end{proof}

\begin{lem}\label{lem:st_case}
    Let $\phi_{+}, \phi_{-}: \Cla \to \Clb(\CH)$ be the unital Jordan $*$-homomorphisms defined by \eqref{equ:bin_Jordan_I}, and $S \in \Clb(\CH)$ be a projection commuting with the ranges of $\alpha$, $\beta$. 
    If there exists a central projection $Q \in \Cla$ such that $S\phi_{+}(Q) = S$ and $S\phi_{-}(Q) = 0$, then there exists a positive contraction $T \in S\Clb(\CH)S$ commuting with the ranges of $\alpha$ and $\beta$ such that 
    \begin{align}\label{equ:st_case_0}
        \begin{pmatrix}
            S\alpha(A) & S\beta(A) \\
            S\beta(A) & S\gamma(A)
        \end{pmatrix} = 
        \begin{pmatrix}
            \sqrt{T} & \sqrt{S-T}\\
            \sqrt{S-T} & -\sqrt{T}    
        \end{pmatrix}
        \begin{pmatrix}
            \phi_{+}(A) & \\
             & \phi_{-}(A)    
        \end{pmatrix}
        \begin{pmatrix}
            \sqrt{T} & \sqrt{S-T}\\
            \sqrt{S-T} & -\sqrt{T}    
        \end{pmatrix}
    \end{align}
    for every $A \in \Cla$. 
\end{lem}

\begin{proof}
Note that $S$ also commutes with the ranges of $\phi_{+}$, $\phi_{-}$, and $\gamma$ (see \eqref{equ:bin_Jordan_I}). 
By replacing $\alpha$, $\beta$, $\gamma$, $\phi_{+}$, and $\phi_{-}$ with $S\alpha$, $S\beta$, $S\gamma$, $S\phi_{+}$, and $S\phi_{-}$, respectively, we may assume that $S = I$. In particular, $\phi_{+}(Q) = I$ and $\phi_{-}(Q) = 0$.

By Lemma \ref{lem:Jordan-completion}, we have $\alpha(A) + \gamma(A) = \phi_{+}(A) + \phi_{-}(A)$ for all $A \in \Cla$. Note that $Q$ is a central projection such that $\phi_{+}(Q) + \phi_{-}(Q) = I$. By Theorem \ref{thm:two_conj_part} and Lemma \ref{lem:bin_Jordan}, we have 
\begin{align*}
\Pi(Q) = \begin{pmatrix}
   T & Y\\
   Y & I-T 
\end{pmatrix},
\end{align*}
where $\Pi$ is the canonical Jordan $*$-homomorphism associated with $\varphi$, $T=\alpha(Q)$ and $Y=\beta(Q)$ are self-adjoint operators that commute with the ranges of $\alpha$, $\beta$, $\phi_{+}$, and $\phi_{-}$. Moreover, since $V_Q = \phi_{+}(Q) - \phi_{-}(Q) \ge 0$, we have $Y \ge 0$. A straightforward computation shows that
\begin{align*}
    W \Pi(Q) W = 
    \begin{pmatrix} 
        I & 0 \\
        0 & 0 
    \end{pmatrix},
\end{align*}
where 
\begin{align*}
    W =  
    \begin{pmatrix} 
        \sqrt{T} & \sqrt{I-T} \\
        \sqrt{I-T} & -\sqrt{T} 
    \end{pmatrix}.
\end{align*}
Consequently, $W \Pi(A) W$ is diagonal. The lemma then follows by applying Lemma \ref{lem:diag_unique}.
\end{proof}

\begin{lem}\label{lem:alpha_normal_Pi_normal}
If $\alpha$ is normal, then the canonical Jordan $*$-homomorphism $\Pi$ associated with the $\Rb$-linear unital unitary preserver $\varphi$ and the Jordan $*$-homomorphisms $\phi_{+}$ and $\phi_{-}$ defined by \eqref{equ:bin_Jordan_I} are normal.
\end{lem}

\begin{proof}
Assume the canonical Jordan $*$-homomorphism $\Pi: \Cla \to \Clb(\CH) \otimes M_2(\Cb)$ admits the form \eqref{equ:Jordan-completion_I}. 
Let $(P_\lambda)$ be a decreasing net of projections in $\Cla$ such that $P_\lambda \downarrow 0$. 
Because $\alpha(P_{\lambda})$ converges to $0$ ultraweakly, the decreasing net of projections $\Pi(P_\lambda)$ must converge to a projection of the form $\left(\begin{smallmatrix} 0 & 0 \\ 0 & F \end{smallmatrix}\right)$. Since
\begin{align*}
    \begin{pmatrix}
        0 & 0 \\
        0 & F
    \end{pmatrix} 
    &\leq \Pi(P_\lambda),
\end{align*}
we have $F \beta(P_\lambda) = 0$ and $F \gamma(P_{\lambda}) = F$ for every $\lambda$. 
Therefore, $FS_{P_{\lambda}}  = 0$, which yields
\begin{align*}
    F = \lim_{\lambda} F \gamma(P_{\lambda}) = \lim_{\lambda} F \left( \alpha(P_\lambda) + S_{P_\lambda}(I-2\alpha(P_\lambda)) \right) = \lim_{\lambda} F \alpha(P_\lambda) = 0.
\end{align*}
Thus, $\Pi$ is normal. Finally, we conclude that both $\phi_{+}$ and $\phi_{-}$ are normal, since 
\begin{align*}
   0 \leq \phi_{\pm}(P_\lambda) \leq \phi_{+}(P_\lambda) + \phi_{-}(P_\lambda) = \alpha(P_{\lambda}) + \gamma(P_\lambda). 
\end{align*}
\end{proof}

Having assembled all the necessary ingredients, we are now ready to prove the main result of this section.

\begin{thm}\label{thm:cls_binary}
Let $\Clc$ be a unital C$^*$-algebra that has no irreducible representation on $\Cb^2$, and $\varphi: \Clc \to \Clb(\CH)$ a purely non-Jordan $\Rb$-linear unital unitary preserver with decomposition
\begin{align*}
    \varphi(H+iK) = \alpha(H)+ i\beta(H) + iJ(\alpha(K) + i\beta(K)), 
\end{align*}
as in Theorem \ref{thm:two_conj_part}. Then there exist unique unital Jordan $*$-homomorphisms $\phi_1, \phi_2 : \Clc \to \Clw^*(\alpha, \beta)$ with mutually commuting ranges and a unique positive contraction $T$ commuting with the ranges of $\alpha$ and $\beta$ such that $JTJ = I-T$, $\supp T(I-T) = I$, and
\begin{align*}
    \alpha(H) = T \phi_{1}(H) + (I-T) \phi_{2}(H), \qquad \beta(H) = \sqrt{T(I-T)}\,(\phi_{1}(H) - \phi_{2}(H))
\end{align*}
for every $H \in \Clc_{\sa}$.
\end{thm}

\begin{proof}
    By the discussion preceding Example \ref{exam:binary_exp} and upon replacing $\Clc$ with $\Cla := \Clc^{**}$, it suffices to establish the existence and uniqueness of $T$, $\phi_1$, and $\phi_2$ under the assumptions that $\Cla$ is a von Neumann algebra without type~$\mathrm{I}_2$ direct summands and that $\varphi: \Cla \to \Clb(\CH)$ is a normal, purely non-Jordan, unital, $\Rb$-linear unitary preserver.

    Let 
    \begin{align*}
        \Cla = \Cla_{\mathrm{I}_1} \oplus \left ( \bigoplus_{n \in \{3, 4, \dots \}} \Cla_{\mathrm{I}_n} \right) \oplus \Cla_{\mathrm{I}_{\infty}} \oplus \Cla_{\mathrm{II}_1} \oplus \Cla_{\mathrm{II}_{\infty}} \oplus \Cla_{\mathrm{III}}
    \end{align*}
    be the decomposition of $\Cla$. For each $s \in \{\mathrm{I}_1, \mathrm{I}_{\infty}, \mathrm{II}_1, \mathrm{II}_{\infty}, \mathrm{III}\} \cup \{\mathrm{I}_n : n \ge 3\}$, let $Q_s$ denote its central projection. By Lemma \ref{lem:alpha_normal_Pi_normal}, the Jordan $*$-homomorphisms $\phi_{+}, \phi_{-}: \Cla \to \Clw^*(\alpha, \beta)$ defined in \eqref{equ:bin_Jordan_I} are normal. For $s, t \in \{\mathrm{I}_1, \mathrm{I}_{\infty}, \mathrm{II}_1, \mathrm{II}_{\infty}, \mathrm{III}\} \cup \{\mathrm{I}_n : n \ge 3\}$, let  
    \begin{align*}
        S_{s,t} \coloneqq \phi_{+}(Q_s) \phi_{-}(Q_t).
    \end{align*}
    Then $\{S_{s, t}\}$ is a family of pairwise orthogonal projections that commute with the ranges of $\alpha$ and $\beta$, and satisfy
    \begin{align*}
        \sum_{s} S_{s,s} + \sum_{s \neq t} S_{s, t} = I.
    \end{align*}
    
    For $s \neq t$, note that $S_{s,t}\phi_{+}(Q_s) = S_{s,t}$ and $S_{s,t}\phi_{-}(Q_s) = 0$.
    By Lemma \ref{lem:st_case}, there exists a positive contraction $T_{s,t} \in S_{s,t}\Clb(\CH)S_{s,t}$ commuting with the ranges of both $\alpha$ and $\beta$ such that
    \begin{align*}
    \begin{pmatrix}
        S_{s,t}\alpha(A) & S_{s,t}\beta(A) \\
        S_{s,t}\beta(A) & S_{s,t}\gamma(A)
    \end{pmatrix} = W_{s,t}
    \begin{pmatrix}
        \phi_{+}(A) & \\
            & \phi_{-}(A)    
    \end{pmatrix}W_{s,t},\qquad \forall A \in \Cla,
    \end{align*}
    where
    \begin{align*}
      W_{s,t} = \begin{pmatrix}
        \sqrt{T_{s,t}} & \sqrt{S_{s,t}-T_{s,t}}\\
        \sqrt{S_{s,t}-T_{s,t}} & -\sqrt{T_{s,t}}    
    \end{pmatrix}.
    \end{align*}
    For a fixed $s$, we have $J S_{s,s} J = S_{s,s}$ by Lemma \ref{lem:bin_Jordan}. Since 
    \begin{align*}
        S_{s,s}\alpha(I-Q_s) + S_{s,s}\gamma(I-Q_s) = S_{s,s}\phi_{+}(I-Q_s) + S_{s,s}\phi_{-}(I-Q_s) = 0,
    \end{align*}
    we obtain $S_{s,s}\alpha(A) = S_{s,s}\beta(A) = 0$ for all $A \in (I-Q_s)\Cla$. In particular, $S_{s,s}\varphi(Q_s) = S_{s,s}$.  And $S_{s,s}\varphi|_{Q_s\Cla}: Q_s\Cla \to \Clb(S_{s,s} \CH)$ is a purely non-Jordan unital $\Rb$-linear unitary preserver. By Lemmas \ref{lem:I1_case} and \ref{lem:case_CM}, there exists a positive contraction $T_{s,s} \in S_{s,s}\Clb(\CH)S_{s,s}$ commuting with the ranges of both $\alpha$ and $\beta$ such that
    \begin{align*}
    \begin{pmatrix}
        S_{s,s}\alpha(A) & S_{s,s}\beta(A) \\
        S_{s,s}\beta(A) & S_{s,s}\gamma(A)
    \end{pmatrix} = W_{s,s}
    \begin{pmatrix}
        \phi_{+}(A) & \\
            & \phi_{-}(A)    
    \end{pmatrix}W_{s,s},\qquad \forall A \in \Cla, 
    \end{align*}
    where 
    \begin{align*}
      W_{s,s} = \begin{pmatrix}
        \sqrt{T_{s,s}} & \sqrt{S_{s,s}-T_{s,s}}\\
        \sqrt{S_{s,s}-T_{s,s}} & -\sqrt{T_{s,s}}    
    \end{pmatrix}.
    \end{align*}
    Thus, by Lemma \ref{lem:diag_unique} the positive contraction    
    \begin{align*}
        T = \sum_{s}T_{s,s} + \sum_{s \neq t}T_{s,t}
    \end{align*} 
    and the Jordan $*$-homomorphisms $\phi_{+}$ and $\phi_{-}$ satisfy the conditions of the theorem.    

    For uniqueness, assume that a positive contraction $T$ and Jordan $*$-homomorphisms $\phi_1, \phi_2$ satisfy the conditions of the theorem. Since $\phi_{1}(P)\phi_2(P) = \phi_2(P)\phi_{1}(P)$ and $\operatorname{supp}\bigl(T(I-T)\bigr) = I$, for every projection $P \in \Cla$ we obtain
    \begin{align*}
        \phi_{+}(P) &= \bigl(T \phi_1(P) + (I-T) \phi_2(P) \bigr)(I-S_P) + \frac{1}{2}(S_P + V_P) = \phi_1(P),\\
        \phi_{-}(P) &= \bigl(T \phi_1(P) + (I-T) \phi_2(P) \bigr)(I-S_P) + \frac{1}{2}(S_P - V_P) = \phi_2(P).
    \end{align*}   
    The uniqueness then follows immediately from Lemma \ref{lem:diag_unique} by noting that 
    \begin{align*}
    \Pi(H) &= 
    \begin{pmatrix}
        T \phi_{1}(H) + (I-T) \phi_{2}(H) & \sqrt{T(I-T)}\,(\phi_{1}(H) - \phi_{2}(H)) \\
        \sqrt{T(I-T)}\,(\phi_{1}(H) - \phi_{2}(H)) & (I-T) \phi_{1}(H) + T \phi_{2}(H)
    \end{pmatrix} \\
    &= \begin{pmatrix}
        \sqrt{T} & \sqrt{I-T} \\
        \sqrt{I-T} & -\sqrt{T}
    \end{pmatrix} 
    \begin{pmatrix}
        \phi_{1}(H) & 0 \\
        0 & \phi_{2}(H)
    \end{pmatrix}
    \begin{pmatrix}
        \sqrt{T} & \sqrt{I-T} \\
        \sqrt{I-T} & -\sqrt{T}
    \end{pmatrix}.
\end{align*}
\end{proof}

\appendix

\section{Proof of Theorem \ref{thm:decmp_general}}

We first establish the following auxiliary fact, which can be used to reduce the proof of Theorem \ref{thm:decmp_general} to more transparent special cases.

\begin{lem}\label{lem:ext_HK_com}
Let $\Clc$ be a unital C$^*$-algebra. Assume that $\rho_1, \rho_2 : \Clc \to \Clb(\CH)$ are norm continuous $\Cb$-linear maps such that 
\begin{align*}
   [\rho_1(H), \rho_2(K)] = 0
\end{align*}    
for every pair of commuting operators $H, K \in \Clc_{\sa}$. Let $\widehat{\rho}_i : \Clc^{**} \to \Clb(\CH)$ be the unique normal extension of $\rho_i$, $i = 1, 2$. Then 
\begin{align*}
   [\widehat{\rho}_1(A), \widehat{\rho}_2(B)] = 0
\end{align*}    
for every pair of commuting operators $A, B \in \Clc^{**}$. 
\end{lem}

\begin{proof}
Let $\omega \in \Clb(\CH)_*$ be a normal functional. Since 
\begin{align*}
    \omega([\rho_1(H), \rho_2(K)]) = 0
\end{align*}
for all $H, K \in \Clc_{\sa}$ satisfying $HK = 0$, Theorem 1.10 in \cite{Goldstein93} (or Theorem 3.1 in \cite{Haag98}) implies there exist $f, g \in \Clc^*$ such that 
\begin{align*}
    \omega([\rho_1(A_1), \rho_2(A_2)]) = f(A_1A_2) + g(A_2A_1)
\end{align*}
for all $A_1, A_2 \in \Clc$. Note that for any $H \in (\Clc)_{\sa}$, we have
\begin{align*}
    f(H) + g(H) = \omega([\rho_1(I), \rho_2(H)]) = 0.
\end{align*}
This implies $g = -f$, and consequently, we have
\begin{align*}
    \omega([\rho_1(A_1), \rho_2(A_2)]) = f([A_1, A_2]).
\end{align*} 
Treating $f$ as its normal extension to $\Clc^{**}$, we obtain
\begin{align*}
    \omega([\widehat{\rho}_1(A), \widehat{\rho}_2(B)]) = f([A, B]) = 0 
\end{align*} 
for every pair of commuting operators $A, B \in \Clc^{**}$. Since this equality holds for every $\omega \in \Clb(\CH)_*$, we conclude that $[\widehat{\rho}_1(A), \widehat{\rho}_2(B)] = 0$. 
\end{proof}

Let $\Cla = \Clc^{**}$. By replacing each $\phi_i$ with its unique normal extension $\widehat{\phi}_i: \Cla \to \Clb(\CH)$, we see from Lemma \ref{lem:ext_HK_com} that it suffices to prove Theorem \ref{thm:decmp_general} when $\Cla$ is a von Neumann algebra with no type I$_2$ direct summand, and $\phi_1, \phi_2: \Cla \to \Clb(\CH)$ are normal unital Jordan $*$-homomorphisms satisfying
\begin{align}\label{equ:HK_commute}
    [\phi_1(S), \phi_2(T)] = 0   
\end{align} 
for all commuting operators $S, T \in \Cla$. Based on this condition, it immediately follows
\begin{align*}
    [\phi_1(Z),\phi_2(A)] = [\phi_1(A),\phi_2(Z)] = 0
\end{align*}
for all $A \in \Cla$ and $Z$ in the center of $\Cla$. Furthermore, by polarization, we have
\begin{align}\label{equ:anti_AB}
    [\phi_1(A), \phi_2(B)] + [\phi_1(B), \phi_2(A)] = 0, \qquad \forall A, B \in \Cla.
\end{align}

Assume that $\Cla$ admits the type decomposition
\begin{align*}
    \Cla = \Cla_{\mathrm{I}_1} \oplus \left ( \bigoplus_{n \in \{3, 4, \dots \}} \Cla_{\mathrm{I}_n} \right) \oplus \Cla_{\mathrm{I}_{\infty}} \oplus \Cla_{\mathrm{II}_1} \oplus \Cla_{\mathrm{II}_{\infty}} \oplus \Cla_{\mathrm{III}}.
\end{align*}
For each $s \in \{\mathrm{I}_1, \mathrm{I}_{\infty}, \mathrm{II}_1, \mathrm{II}_{\infty}, \mathrm{III}\} \cup \{\mathrm{I}_n : n \ge 3\}$, we identify $\Cla_s$ with its corresponding direct summand in $\Cla$, and let $Q_s$ denote its central projection. We set $Q_s = 0$ if $\Cla$ has no direct summand of type $s$.

We claim that Theorem \ref{thm:decmp_general} is an immediate corollary of the following special case.

\begin{thm}\label{thm:one_sum_case}
Let $\Clm = \Clc \otimes M_n(\Cb)$, where $\Clc$ is a unital C$^*$-algebra and $n \geq 3$. Let $\phi_1, \phi_2 : \Clm \to \Clb(\CH)$ be unital Jordan $*$-homomorphisms such that  
\begin{align*}
    [\phi_1(A), \phi_2(B)] = 0
\end{align*}
for all commuting operators $A, B \in \Clm$. Then there exists a projection $P \in \Clb(\CH)$ commuting with both ranges $\phi_1(\Clm)$ and $\phi_2(\Clm)$ such that
\begin{align*}
    \phi_1(A_1)P = \phi_2(A_1)P \quad \text{and} \quad [(I-P)\phi_1(A_1), (I-P)\phi_2(A_2)] = 0
\end{align*}
for all $A_1, A_2 \in \Clm$. 
\end{thm}

Indeed, using the notation introduced at the beginning of the appendix, the normality of $\phi_1$ and $\phi_2$ implies that
\begin{align*}
    \sum_s \phi_i(Q_s) = I, \qquad i = 1, 2.
\end{align*}
Set
\begin{align*}
    P_0 := \sum_{s \neq \mathrm{I}_1} \phi_1(Q_s) \phi_2(Q_s). 
\end{align*}
Since $[\phi_1(L_1), \phi_2(L_2)] = 0$ for all $L_1, L_2 \in \Cla_{\mathrm{I}_1}$, and   
\begin{align*}
    [\phi_1(A_s), \phi_2(B_t)] = 0 \qquad \text{for all } A_s \in \Cla_s, \, B_t \in \Cla_t \text{ with } s \neq t,
\end{align*}
it follows that
\begin{align*}
    [(I-P_0)\phi_1(A), (I-P_0)\phi_2(B)] = 0 \qquad \text{for all } A, B \in \Cla.
\end{align*}
For each $s \neq \mathrm{I}_1$, applying Theorem \ref{thm:one_sum_case} to the restricted maps on $\Cla_s$ yields a projection $P_s \leq \phi_1(Q_s)\phi_2(Q_s)$ commuting with both $\phi_1(\Cla)$ and $\phi_2(\Cla)$ such that $\phi_1(A_1)P_s = \phi_2(A_1)P_s$ and 
\begin{align*}
    \bigl[ (\phi_1(Q_s)\phi_2(Q_s) - P_s)\phi_1(A_1), \, (\phi_1(Q_s)\phi_2(Q_s) - P_s)\phi_2(A_2) \bigr] = 0
\end{align*}
for all $A_1, A_2 \in \Cla$. Consequently, the projection
\begin{align*}
    P := \sum_{s \neq \mathrm{I}_1} P_s 
\end{align*}
satisfies the conditions of Theorem \ref{thm:decmp_general}. The remainder of this appendix is devoted to the proof of Theorem \ref{thm:one_sum_case}. While normality of $\phi_1$ and $\phi_2$ is used to reduce Theorem \ref{thm:decmp_general} to Theorem \ref{thm:one_sum_case}, Theorem \ref{thm:one_sum_case} itself holds without this assumption.

\begin{lem}\label{lem:cen_decmp_comp_gen}
Let $\phi_1, \phi_2 : \Clm = \Clc \otimes M_n(\Cb) \to \Clb(\CH)$ be unital Jordan $*$-homomorphisms satisfying the hypotheses of Theorem \ref{thm:one_sum_case}. If $P$ is a central projection in the von Neumann algebra generated by the image of $\phi_1$ such that $P\phi_1$ is a $*$-homomorphism and $(I-P)\phi_1$ is a  $*$-anti-homomorphism, then 
\begin{align*}
    [P, \phi_2(A)] = 0, \qquad \forall A \in \Clm.
\end{align*}
\end{lem}

\begin{proof}
Define the map $\rho$ on $\Clm$ by
\begin{align*}
    \rho(T) := P\phi_2(T)(I-P), \quad T \in \Clm.
\end{align*}
It suffices to show that $\rho$ vanishes identically on $\Clm$. To this end, we first prove that $\rho(I \otimes E) = 0$ for every rank-one projection $E \in M_n(\Cb)$. For every $A, B \in \bigl ( I \otimes (I-E) \bigr )\Clm \bigl (I \otimes (I-E) \bigr )$, we have
\begin{align*}
   \rho(I \otimes E)\phi_1(AB) &= \phi_1(AB) \rho(I \otimes E) = \phi_1(A)\phi_1(B)\rho(I \otimes E)\\
   & = \rho(I \otimes E) (I-P) \phi_1(A)\phi_1(B) = \rho(I \otimes E)\phi_1(BA). 
\end{align*} 
Therefore, $\rho(I \otimes E) \phi_1([A, B]) = 0$ for every $A, B \in \bigl ( I \otimes (I-E) \bigr )\Clm \bigl (I \otimes (I-E) \bigr )$. Furthermore, for any $C, D \in \bigl ( I \otimes (I-E) \bigr )\Clm \bigl (I \otimes (I-E) \bigr )$, this implies
\begin{equation}\label{equ:cen_decmp_comp_gen_0}
   \begin{aligned}
     \rho(I \otimes E) \phi_1(C[A, B]D) &= \rho(I \otimes E) \phi_1(D)\phi_1([A, B])\phi_1(C)\\ 
    &= \phi_1(D)  \rho(I \otimes E)\phi_1([A, B])\phi_1(C) = 0.
   \end{aligned} 
\end{equation}
Note that the closed ideal generated by the commutators in $(I-E)M_n(\Cb)(I-E) (\simeq M_{n-1}(\Cb))$ is the entire algebra. We have
\begin{align}\label{equ:cen_decmp_comp_gen_I}
   \rho(I \otimes E) = \rho(I \otimes E) \phi_1(I \otimes E).
\end{align}
Let $\{E_1, E_2, \ldots, E_n\} \subset M_n(\Cb)$ be a maximal family of mutually orthogonal rank-one projections with $E_1 = E$. By Equation \eqref{equ:cen_decmp_comp_gen_I}, 
\begin{align*}
   0 = \rho(I) = \sum_{k=1}^n \rho(I \otimes E_k) \phi_1(I \otimes E_k).
\end{align*}
Then 
\begin{align*}
  \rho(I \otimes E ) = \rho(I \otimes E) \phi_1(I \otimes E) =  \left ( \sum_{k=1}^n \rho(I \otimes E_k) \phi_1(I \otimes E_k) \right ) \phi_1(I \otimes E) = 0.
\end{align*}
Since every operator in $M_n(\Cb)$ is a linear combination of rank-one projections, we have
\begin{align}\label{equ:cen_decmp_comp_gen_II}
    \rho(S) = 0, \quad \forall S \in I \otimes M_n(\Cb).
\end{align}

Let $T \in \Clm$ and $S \in I \otimes M_n(\Cb)$. By \eqref{equ:cen_decmp_comp_gen_II}, we obtain
\begin{align*}
  0 = \phi_1(T + S)\rho(T + S) - \rho(T + S)\phi_1(T + S) = \phi_1(S)\rho(T) - \rho(T)\phi_1(S).  
\end{align*}
Then the same computation used to obtain \eqref{equ:cen_decmp_comp_gen_0} shows that
\begin{align*}
    \rho(T)\phi_1(X) = 0
\end{align*}
for every $X$ in the ideal of $I \otimes M_n(\Cb)$ generated by its commutators. In particular, $\rho(T) = 0$.
\end{proof}

\begin{lem}\label{lem:In_case}
    Let $\rho: M_n(\Cb) \to \Clb(\CK) \otimes M_n(\Cb)$ be a unital Jordan $*$-homomorphism, where $\CK$ is a Hilbert space and $n \geq 3$. Assume that 
   \begin{align*}
         [I \otimes A, \rho(B)] = 0 
   \end{align*}   
   whenever $A, B \in M_n(\Cb)$ commute. Then there is a projection $Q \in \Clb(\CK)$ and a unital Jordan $*$-homomorphism $\pi: M_n(\Cb) \to (I-Q)\Clb(\CK)(I-Q)$ such that  
    \begin{align*}
        \rho(T)= Q \otimes T + \pi(T) \otimes I_n, \qquad \forall T \in M_n(\Cb).
    \end{align*}
\end{lem}

\begin{proof}
For every rank-one projection $E \in M_n(\Cb)$, $\rho(E)$ is in the commutant of $I \otimes (I-E)M_n(\Cb)(I-E)$. Therefore, 
\begin{align*}
   \rho(E) = Q(E) \otimes E + \pi(E) \otimes I_n, 
\end{align*}
where $Q(E), \pi(E) \in (\Clb(\CK))_{\sa}$. 

Let $\{E_1, E_2, \ldots, E_n\} \subset M_n(\Cb)$ be a family of mutually orthogonal rank-one projections. For each $k \in \{1, \ldots, n\}$, 
\begin{align*}
    I \otimes E_k &= \left( \sum_{j=1}^n Q(E_j) \otimes E_j + \pi(E_j) \otimes I_n \right) (I \otimes E_k) = \left( Q(E_k) + \sum_{j=1}^n \pi(E_j) \right) \otimes E_k.
\end{align*}
Therefore, $Q(E_k) = I - \sum_{j=1}^n \pi(E_j) = Q(E_j)$ for all $j, k \in \{1, \ldots, n\}$. Given any two rank-one projections $E, F \in M_n(\Cb)$, there exists a rank-one projection $P$ orthogonal to both. Then the above argument implies that $Q(E) = Q(P) = Q(F)$. We may therefore write simply $Q$ for $Q(E)$.

Note that $\{\rho(E_j)\}_{j=1}^n$ is a family of mutually orthogonal projections whenever $\{E_1, E_2, \ldots, E_n\}$ is a family of mutually orthogonal rank-one projections. The identity $\rho(E_j) \rho(E_k) = \delta_{j,k}\rho(E_k)$ then yields
\begin{align*}
    \pi(E_j) \pi(E_k) = \delta_{j,k}\pi(E_k),
\end{align*}
as well as $(Q + \pi(E_j)) \pi(E_k) = 0$ for all $j \neq k$. Therefore, each $\pi(E_j)$ is a projection and $Q\pi(E_j) = 0$ for $j = 1, \ldots, n$. Furthermore, the identity $\rho(E_j)^2 = \rho(E_j)$ implies that $Q$ is also a projection. In summary, we have
\begin{align*}
    \rho(E) = Q \otimes E + \pi(E) \otimes I_n,
\end{align*}
where $Q$ and $\pi(E)$ are mutually orthogonal projections in $\Clb(\CK)$. The lemma now follows by linearity, since $M_n(\Cb)$ is linearly spanned by rank-one projections. And the unitality of $\pi$ is ensured by that of $\rho$.
\end{proof}

\begin{proof}[\textbf{Proof of Theorem \ref{thm:one_sum_case}}]
Let $F$ be the central projection in the von Neumann algebra generated by $\phi_1(\Clm)$ such that $F \phi_1$ is a $*$-homomorphism and $(I - F)\phi_1$ is a $*$-anti-homomorphism. By Lemma~\ref{lem:cen_decmp_comp_gen}, both $F \phi_2$ and $(I-F)\phi_2$ are unital Jordan $*$-homomorphisms. By considering the pairs $\{F \phi_1, F\phi_2\}$ and $\{(I-F)\phi_1, (I-F)\phi_2\}$ separately, and replacing $\Clm$ by $\Clm^{\op}$ if necessary, we may assume without loss of generality that $\phi_1$ is a unital $*$-homomorphism.

Let $\{E_{ij}\}_{1 \leq i, j\leq n}$ be the standard matrix units of $M_n(\Cb)$. Without loss of generality, we may assume that $\Clb(\CH) = \Clb(\CK) \otimes M_n(\Cb)$, and  
\begin{align*}
    \phi_1(T \otimes E_{ij}) = \sigma_1(T) \otimes E_{ij},
\end{align*}
where $\sigma_1: \Clc \to \Clb(\CK)$ is a $*$-homomorphism. Applying Lemma \ref{lem:In_case} to $\phi_2|_{I \otimes M_n(\Cb)}$ yields a projection $Q \in \Clb(\CK)$ and a unital Jordan $*$-homomorphism $\pi: M_n(\Cb) \to (I-Q) \Clb(\CK) (I-Q)$ such that
\begin{align*}
    \phi_2(I \otimes B) = Q \otimes B + \pi(B) \otimes I_n, \qquad B \in M_n(\Cb).
\end{align*}  
Since $T \otimes I$ commutes with $I \otimes B$ for every $T \in \Clc$ and $B \in M_n(\Cb)$, there exists a Jordan $*$-homomorphism $\sigma_2: \Clc \to \Clb(\CK)$ such that 
\begin{align*}
    \phi_2(T \otimes I) = \sigma_2(T) \otimes I_n.
\end{align*}

For every rank-one projection $E \in M_n(\Cb)$, the commutativity of $\phi_2(I \otimes E)$ and $\phi_s(T \otimes I)$ ($s = 1, 2$) implies that
\begin{align*}
   (Q + \pi(E)) \sigma_s(T) = \sigma_s(T) (Q + \pi(E))\quad \text{and} \quad \pi(E)\sigma_s(T) = \sigma_s(T)\pi(E). 
\end{align*}
Therefore, we obtain 
\begin{align*}
Q \sigma_s(T) = \sigma_s(T)Q, \qquad \pi(B) \sigma_s(T) = \sigma_s(T) \pi(B)
\end{align*}
for every $T \in \Clc$ and $B \in M_n(\Cb)$.

Let $P = Q \otimes I$. Then $P$ commutes with the ranges of $\phi_1$ and $\phi_2$. And the identity
\begin{align*}
    0 = P[\phi_1(T \otimes E_{12} + I \otimes E_{23}), \phi_2(T \otimes E_{12} + I \otimes E_{23})]P = (\sigma_1(T) - \sigma_2(T)) Q \otimes E_{13}
\end{align*}
implies that $\sigma_1(T)Q = \sigma_2(T)Q$, which in turn implies that $P\phi_1 = P\phi_2$. Finally, by Equation \eqref{equ:anti_AB}, we have 
\begin{align*}
   [\sigma_1(T_1), \sigma_2(T_2)]\pi(E_{23}) \otimes E_{12} &=  [(I-P)\phi_1(T_1 \otimes E_{12}), (I-P)\phi_2(T_2 \otimes E_{23})]\\
   &= -[(I-P)\phi_1(T_2 \otimes E_{23}), (I-P)\phi_2(T_1 \otimes E_{12})]\\
   &= [\sigma_2(T_1), \sigma_1(T_2)]\pi(E_{12}) \otimes E_{23},
\end{align*}
which implies that $[(I-Q)\sigma_1(T_1), (I-Q)\sigma_2(T_2)] = 0$. It follows that 
\begin{align*}
    [(I-P)\phi_1(A_1), (I-P)\phi_2(A_2)] = 0, \qquad \forall A_1, A_2 \in \Clm.
\end{align*}
\end{proof}

\printbibliography

@article{Sta82,
    author = {P. J. Stacey},
    title = {Type {$I_2$} {JBW}-algebras},
    journal = {Quart. J. Math.},
    volume = {33},
    year = {1982},
    number = {1},
    pages = {115--127},
 }

@article{kl57,
  author  = {D. C. Kleinecke},
  title   = {On operator commutators},
  journal = {Proc. Amer. Math. Soc.},
  volume  = {8},
  pages   = {535--536},
  year    = {1957},
}

@article{ShiShenGuMa2022,
  author  = {W. Shi and J. Shen and W. Gu and M. Ma},
  title   = {{$L^p$}-isometries of Grassmann spaces in factors of type {II}},
  journal = {J. Operator Theory},
  volume  = {87},
  year    = {2022},
  pages   = {389--412},
}

@article{QianShenShiWuYuan2023,
  author  = {W. Qian and J. Shen and W. Shi and W. Wu and W. Yuan},
  title   = {Surjective {$L^p$}-isometries on Grassmann spaces},
  journal = {Sci. China Math.},
  volume  = {66},
  year    = {2023},
  pages   = {2105--2118},
}

@article{GuWuYuan2020,
  author  = {W. Gu and W. Wu and W. Yuan},
  title   = {Transition probability preserving maps on a Grassmann space in a semifinite factor},
  journal = {J. Math. Anal. Appl.},
  volume  = {487},
  year    = {2020},
  pages   = {123957}
}

@article{QianWangWuYuan2019,
  author  = {W. Qian and L. Wang and W. Wu and W. Yuan},
  title   = {{Wigner}-type theorem on transition probability preserving maps in semifinite factors},
  journal = {J. Funct. Anal.},
  volume  = {276},
  year    = {2019},
  pages   = {1773--1787}
}

@article{Mori2019,
  author  = {M. Mori},
  title   = {Isometries between projection lattices of von {Neumann} algebras},
  journal = {J. Funct. Anal.},
  volume  = {276},
  year    = {2019},
  pages   = {3511--3528}
}

@article{GeherSemrl2016,
  author  = {G. P. Geh{\'e}r and P. {\v{S}}emrl},
  title   = {Isometries of Grassmann spaces},
  journal = {J. Funct. Anal.},
  volume  = {270},
  year    = {2016},
  pages   = {1585--1601}
}

@article{GeherSemrl2018,
  author  = {G. P. Geh{\'e}r and P. {\v{S}}emrl},
  title   = {Isometries of Grassmann spaces, {II}},
  journal = {Adv. Math.},
  volume  = {332},
  year    = {2018},
  pages   = {287--310}
}

@article{Dye1955,
  author  = {H. A. Dye},
  title   = {On the geometry of projections in certain operator algebras},
  journal = {Ann. of Math. (2)},
  volume  = {61},
  year    = {1955},
  pages   = {73--89}
}

@article{LiPierce2001,
  author  = {C. K. Li and S. Pierce},
  title   = {Linear preserver problems},
  journal = {Amer. Math. Monthly},
  volume  = {108},
  number = {7},
  year    = {2001},
  pages   = {591--605}
}

@book{BM04,
  author    = {D. P. Blecher and C. L. Merdy},
  title     = {Operator Algebras and Their Modules---An Operator Space Approach},
  series    = {London Mathematical Society Monographs. New Series},
  volume    = {30},
  publisher = {Oxford University Press},
  address   = {Oxford},
  year      = {2004},
  isbn      = {0-19-852659-8}
}

@article{Goldstein93,
  author  = {S. Goldstein},
  title   = {Stationarity of operator algebras},
  journal = {Journal of Functional Analysis},
  year    = {1993},
  volume  = {118},
  number  = {2},
  pages   = {275--308},
}

@incollection{Haag98,
  author    = {U. Haagerup and N. J. Laustsen},
  title     = {Weak amenability of {$C^*$}-algebras and a theorem of Goldstein},
  booktitle = {Banach Algebras '97},
  publisher = {Walter de Gruyter},
  address   = {Berlin},
  year      = {1998},
  pages     = {223--243}
}

@article{bunce92,
  author  = {L. J. Bunce and J. D. Wright},
  title   = {The {Mackey}-{Gleason} Problem},
  journal = {Bulletin of the American Mathematical Society},
  volume  = {26},
  number  = {2},
  pages   = {288--293},
  year    = {1992},
  month   = {4},
}

@article{kadison85,
  title={Means and convex combinations of unitary operators},
  author={R. V. Kadison and G. K. Pedersen},
  journal={Mathematica Scandinavica},
  volume={57},
  pages={249--266},
  year={1985},
}

@article{stormer65,
  author  = {E. St{\o}rmer},
  title   = {On the Jordan structure of {C*}-algebras},
  journal = {Transactions of the American Mathematical Society},
  volume  = {120},
  number  = {3},
  pages   = {438--447},
  year    = {1965},
  publisher = {American Mathematical Society}
}

@article{choi87,
  title={Linear maps preserving commutativity},
  author={M. D. Choi and A. A. Jafarian and H. Radjavi},
  journal={Linear Algebra and its Applications},
  volume={87},
  pages={227--242},
  year={1987},
  publisher={Elsevier}
}

@book{OS84,
  title={Jordan Operator Algebras},
  author={H. Hanche-Olsen and E. St{\o}rmer},
  year={1984},
  publisher={Pitman},
  series={Monographs and Studies in Mathematics},
  volume={21}
}

@article{Wei1975,
  author  = {A. Wei},
  title   = {Linear Transformations Preserving the Real Orthogonal Group},
  journal = {Canadian Journal of Mathematics},
  volume  = {27},
  number  = {3},
  pages   = {561--572},
  year    = {1975},
}

@article{marcus59,
  title={All linear operators leaving the unitary group invariant},
  author={M. Marcus},
  journal={Duke Mathematical Journal},
  volume={26},
  number={1},
  pages={155--163},
  year={1959},
  publisher={Duke University Press}
}

@book{M2007,
  title     = {Selected Preserver Problems on Algebraic Structures of Linear Operators and on Function Spaces},
  author    = {L. Moln{\'a}r},
  series    = {Lecture Notes in Mathematics},
  volume    = {1895},
  year      = {2007},
  publisher = {Springer-Verlag},
  address   = {Berlin, Heidelberg},
  isbn      = {978-3-540-39944-5}
}

@Article{RD66,
  author    = {B. Russo and H.A. Dye},
  journal   = {Duke Math. J.},
  title     = {A note on unitary operators in {C}$^*$-algebras},
  year      = {1966},
  pages     = {413–416},
  volume    = {33},
}

@Book{KR97,
  author    = {R. V. Kadison and J. R. Ringrose},
  publisher = {American Mathematical Society},
  title     = {Fundamentals of the theory of operator algebras, Vol.\ I, II},
  series    ={Graduate Studies in Mathematics},
  volume    = {15, 21},
  year      = {1997},
}

@Book{T02,
  author    = {M. Takesaki},
  publisher = {Springer-Verlag Berlin Heidelberg},
  title     = {Theory of Operator Algebra I},
  series    ={Encyclopaedia of Mathematical Sciences},
  volume    = {124},
  year      = {2002},
}

\end{document}